\input amstex

\input texdraw
\input epsf

\documentstyle{amsppt}
%\magnification=\magstep1
\pagewidth{5.5truein}
\hcorrection{0.55in}
\pageheight{7.5truein}\vcorrection{0.75in}
%\pagewidth{6.48truein}%\hcorrection{0.55in}
%\pageheight{9.0truein}%\vcorrection{0.75in}
\TagsOnRight
\NoRunningHeads

\catcode`\@=11
\font@\twelverm=cmr10 scaled\magstep1
\font@\twelveit=cmti10 scaled\magstep1
\font@\twelvebf=cmbx10 scaled\magstep1
\font@\twelvei=cmmi10 scaled\magstep1
\font@\twelvesy=cmsy10 scaled\magstep1
\font@\twelveex=cmex10 scaled\magstep1

\font@\fourteenrm=cmr10 scaled\magstep2
\font@\fourteenit=cmti10 scaled\magstep2
\font@\fourteensl=cmsl10 scaled\magstep2
\font@\fourteensmc=cmcsc10 scaled\magstep2
\font@\fourteentt=cmtt10 scaled\magstep2
\font@\fourteenbf=cmbx10 scaled\magstep2
\font@\fourteeni=cmmi10 scaled\magstep2
\font@\fourteensy=cmsy10 scaled\magstep2
\font@\fourteenex=cmex10 scaled\magstep2
\font@\fourteenmsa=msam10 scaled\magstep2
\font@\fourteeneufm=eufm10 scaled\magstep2
\font@\fourteenmsb=msbm10 scaled\magstep2
\newtoks\fourteenpoint@
\def\fourteenpoint{\normalbaselineskip15\p@
 \abovedisplayskip18\p@ plus4.3\p@ minus12.9\p@
 \belowdisplayskip\abovedisplayskip
 \abovedisplayshortskip\z@ plus4.3\p@
 \belowdisplayshortskip10.1\p@ plus4.3\p@ minus5.8\p@
 \textonlyfont@\rm\fourteenrm \textonlyfont@\it\fourteenit
 \textonlyfont@\sl\fourteensl \textonlyfont@\bf\fourteenbf
 \textonlyfont@\smc\fourteensmc \textonlyfont@\tt\fourteentt
%Erg\"anzung des fetten Small-Capitals-Fonts:
%
 \ifsyntax@ \def\big##1{{\hbox{$\left##1\right.$}}}%
  \let\Big\big \let\bigg\big \let\Bigg\big
 \else
  \textfont\z@=\fourteenrm  \scriptfont\z@=\twelverm  \scriptscriptfont\z@=\tenrm
  \textfont\@ne=\fourteeni  \scriptfont\@ne=\twelvei  \scriptscriptfont\@ne=\teni
  \textfont\tw@=\fourteensy \scriptfont\tw@=\twelvesy \scriptscriptfont\tw@=\tensy
  \textfont\thr@@=\fourteenex \scriptfont\thr@@=\twelveex
        \scriptscriptfont\thr@@=\twelveex
  \textfont\itfam=\fourteenit \scriptfont\itfam=\twelveit
        \scriptscriptfont\itfam=\twelveit
  \textfont\bffam=\fourteenbf \scriptfont\bffam=\twelvebf
        \scriptscriptfont\bffam=\tenbf
  \setbox\strutbox\hbox{\vrule height12.2\p@ depth5\p@ width\z@}%
  \setbox\strutbox@\hbox{\lower.72\normallineskiplimit\vbox{%
        \kern-\normallineskiplimit\copy\strutbox}}%
 \setbox\z@\vbox{\hbox{$($}\kern\z@}\bigsize@=1.7\ht\z@
 \fi
 \normalbaselines\rm\ex@.2326ex\jot4.3\ex@\the\fourteenpoint@}
\catcode`\@=13

\catcode`\@=11
\font\tenln    = line10
\font\tenlnw   = linew10

\newskip\Einheit \Einheit=0.5cm
\newcount\xcoord \newcount\ycoord
\newdimen\xdim \newdimen\ydim \newdimen\PfadD@cke \newdimen\Pfadd@cke

%%%%%%%%%%%%%%%%%%%%%%%%%%%%%%%%%%%%%%%%%%%%%%%%%
%LaTeX counters, dimensions, variables for lines%
%%%%%%%%%%%%%%%%%%%%%%%%%%%%%%%%%%%%%%%%%%%%%%%%%
\newcount\@tempcnta
\newcount\@tempcntb

\newdimen\@tempdima
\newdimen\@tempdimb

\newdimen\@wholewidth
\newdimen\@halfwidth

\newcount\@xarg
\newcount\@yarg
\newcount\@yyarg
\newbox\@linechar
\newbox\@tempboxa
\newdimen\@linelen
\newdimen\@clnwd
\newdimen\@clnht

\newif\if@negarg

\def\@whilenoop#1{}
\def\@whiledim#1\do #2{\ifdim #1\relax#2\@iwhiledim{#1\relax#2}\fi}
\def\@iwhiledim#1{\ifdim #1\let\@nextwhile=\@iwhiledim
        \else\let\@nextwhile=\@whilenoop\fi\@nextwhile{#1}}

\def\@whileswnoop#1\fi{}
\def\@whilesw#1\fi#2{#1#2\@iwhilesw{#1#2}\fi\fi}
\def\@iwhilesw#1\fi{#1\let\@nextwhile=\@iwhilesw
         \else\let\@nextwhile=\@whileswnoop\fi\@nextwhile{#1}\fi}

\def\thinlines{\let\@linefnt\tenln \let\@circlefnt\tencirc
  \@wholewidth\fontdimen8\tenln \@halfwidth .5\@wholewidth}
\def\thicklines{\let\@linefnt\tenlnw \let\@circlefnt\tencircw
  \@wholewidth\fontdimen8\tenlnw \@halfwidth .5\@wholewidth}
\thinlines
%%%%%%%%%%%%%%%%%%%%%%%%%%%%%%%%%%%%%%%%%%%%%%%%%%%%%%%%%%%

\PfadD@cke1pt \Pfadd@cke0.5pt
\def\PfadDicke#1{\PfadD@cke#1 \divide\PfadD@cke by2 \Pfadd@cke\PfadD@cke \multiply\PfadD@cke by2}
\long\def\LOOP#1\REPEAT{\def\BODY{#1}\ITERATE}
\def\ITERATE{\BODY \let\next\ITERATE \else\let\next\relax\fi \next}
\let\REPEAT=\fi
\def\Punkt{\hbox{\raise-2pt\hbox to0pt{\hss$\ssize\bullet$\hss}}}
\def\DuennPunkt(#1,#2){\unskip
  \raise#2 \Einheit\hbox to0pt{\hskip#1 \Einheit
          \raise-2.5pt\hbox to0pt{\hss$\bullet$\hss}\hss}}
\def\NormalPunkt(#1,#2){\unskip
  \raise#2 \Einheit\hbox to0pt{\hskip#1 \Einheit
          \raise-3pt\hbox to0pt{\hss\twelvepoint$\bullet$\hss}\hss}}
\def\DickPunkt(#1,#2){\unskip
  \raise#2 \Einheit\hbox to0pt{\hskip#1 \Einheit
          \raise-4pt\hbox to0pt{\hss\fourteenpoint$\bullet$\hss}\hss}}
\def\Kreis(#1,#2){\unskip
  \raise#2 \Einheit\hbox to0pt{\hskip#1 \Einheit
          \raise-4pt\hbox to0pt{\hss\fourteenpoint$\circ$\hss}\hss}}

%%%%%%%%%%%%%%%%%%%%%
%LaTeX line macros%
%%%%%%%%%%%%%%%%%%%%%
\def\Line@(#1,#2)#3{\@xarg #1\relax \@yarg #2\relax
\@linelen=#3\Einheit
\ifnum\@xarg =0 \@vline
  \else \ifnum\@yarg =0 \@hline \else \@sline\fi
\fi}

\def\@sline{\ifnum\@xarg< 0 \@negargtrue \@xarg -\@xarg \@yyarg -\@yarg
  \else \@negargfalse \@yyarg \@yarg \fi
\ifnum \@yyarg >0 \@tempcnta\@yyarg \else \@tempcnta -\@yyarg \fi
\ifnum\@tempcnta>6 \@badlinearg\@tempcnta0 \fi
\ifnum\@xarg>6 \@badlinearg\@xarg 1 \fi
\setbox\@linechar\hbox{\@linefnt\@getlinechar(\@xarg,\@yyarg)}%
\ifnum \@yarg >0 \let\@upordown\raise \@clnht\z@
   \else\let\@upordown\lower \@clnht \ht\@linechar\fi
\@clnwd=\wd\@linechar
\if@negarg \hskip -\wd\@linechar \def\@tempa{\hskip -2\wd\@linechar}\else
     \let\@tempa\relax \fi
\@whiledim \@clnwd <\@linelen \do
  {\@upordown\@clnht\copy\@linechar
   \@tempa
   \advance\@clnht \ht\@linechar
   \advance\@clnwd \wd\@linechar}%
\advance\@clnht -\ht\@linechar
\advance\@clnwd -\wd\@linechar
\@tempdima\@linelen\advance\@tempdima -\@clnwd
\@tempdimb\@tempdima\advance\@tempdimb -\wd\@linechar
\if@negarg \hskip -\@tempdimb \else \hskip \@tempdimb \fi
\multiply\@tempdima \@m
\@tempcnta \@tempdima \@tempdima \wd\@linechar \divide\@tempcnta \@tempdima
\@tempdima \ht\@linechar \multiply\@tempdima \@tempcnta
\divide\@tempdima \@m
\advance\@clnht \@tempdima
\ifdim \@linelen <\wd\@linechar
   \hskip \wd\@linechar
  \else\@upordown\@clnht\copy\@linechar\fi}

\def\@hline{\ifnum \@xarg <0 \hskip -\@linelen \fi
\vrule height\Pfadd@cke width \@linelen depth\Pfadd@cke
\ifnum \@xarg <0 \hskip -\@linelen \fi}

\def\@getlinechar(#1,#2){\@tempcnta#1\relax\multiply\@tempcnta 8
\advance\@tempcnta -9 \ifnum #2>0 \advance\@tempcnta #2\relax\else
\advance\@tempcnta -#2\relax\advance\@tempcnta 64 \fi
\char\@tempcnta}

\def\Vektor(#1,#2)#3(#4,#5){\unskip\leavevmode
  \xcoord#4\relax \ycoord#5\relax
      \raise\ycoord \Einheit\hbox to0pt{\hskip\xcoord \Einheit
         \Vector@(#1,#2){#3}\hss}}

\def\Vector@(#1,#2)#3{\@xarg #1\relax \@yarg #2\relax
\@tempcnta \ifnum\@xarg<0 -\@xarg\else\@xarg\fi
\ifnum\@tempcnta<5\relax
\@linelen=#3\Einheit
\ifnum\@xarg =0 \@vvector
  \else \ifnum\@yarg =0 \@hvector \else \@svector\fi
\fi
\else\@badlinearg\fi}

\def\@hvector{\@hline\hbox to 0pt{\@linefnt
\ifnum \@xarg <0 \@getlarrow(1,0)\hss\else
    \hss\@getrarrow(1,0)\fi}}

\def\@vvector{\ifnum \@yarg <0 \@downvector \else \@upvector \fi}

\def\@svector{\@sline
\@tempcnta\@yarg \ifnum\@tempcnta <0 \@tempcnta=-\@tempcnta\fi
\ifnum\@tempcnta <5
  \hskip -\wd\@linechar
  \@upordown\@clnht \hbox{\@linefnt  \if@negarg
  \@getlarrow(\@xarg,\@yyarg) \else \@getrarrow(\@xarg,\@yyarg) \fi}%
\else\@badlinearg\fi}

\def\@upline{\hbox to \z@{\hskip -.5\Pfadd@cke \vrule width \Pfadd@cke
   height \@linelen depth \z@\hss}}

\def\@downline{\hbox to \z@{\hskip -.5\Pfadd@cke \vrule width \Pfadd@cke
   height \z@ depth \@linelen \hss}}

\def\@upvector{\@upline\setbox\@tempboxa\hbox{\@linefnt\char'66}\raise
     \@linelen \hbox to\z@{\lower \ht\@tempboxa\box\@tempboxa\hss}}

\def\@downvector{\@downline\lower \@linelen
      \hbox to \z@{\@linefnt\char'77\hss}}

\def\@getlarrow(#1,#2){\ifnum #2 =\z@ \@tempcnta='33\else
\@tempcnta=#1\relax\multiply\@tempcnta \sixt@@n \advance\@tempcnta
-9 \@tempcntb=#2\relax\multiply\@tempcntb \tw@
\ifnum \@tempcntb >0 \advance\@tempcnta \@tempcntb\relax
\else\advance\@tempcnta -\@tempcntb\advance\@tempcnta 64
\fi\fi\char\@tempcnta}

\def\@getrarrow(#1,#2){\@tempcntb=#2\relax
\ifnum\@tempcntb < 0 \@tempcntb=-\@tempcntb\relax\fi
\ifcase \@tempcntb\relax \@tempcnta='55 \or
\ifnum #1<3 \@tempcnta=#1\relax\multiply\@tempcnta
24 \advance\@tempcnta -6 \else \ifnum #1=3 \@tempcnta=49
\else\@tempcnta=58 \fi\fi\or
\ifnum #1<3 \@tempcnta=#1\relax\multiply\@tempcnta
24 \advance\@tempcnta -3 \else \@tempcnta=51\fi\or
\@tempcnta=#1\relax\multiply\@tempcnta
\sixt@@n \advance\@tempcnta -\tw@ \else
\@tempcnta=#1\relax\multiply\@tempcnta
\sixt@@n \advance\@tempcnta 7 \fi\ifnum #2<0 \advance\@tempcnta 64 \fi
\char\@tempcnta}
%%%%%%%%%%%%%%%%%%%%%%%%%%%%%%%%%%%%%%%%%%%%%%%%%%%%%%%%%%%%%

\def\Diagonale(#1,#2)#3{\unskip\leavevmode
  \xcoord#1\relax \ycoord#2\relax
      \raise\ycoord \Einheit\hbox to0pt{\hskip\xcoord \Einheit
         \Line@(1,1){#3}\hss}}
\def\AntiDiagonale(#1,#2)#3{\unskip\leavevmode
  \xcoord#1\relax \ycoord#2\relax %\advance\xcoord by -0.05\relax
      \raise\ycoord \Einheit\hbox to0pt{\hskip\xcoord \Einheit
         \Line@(1,-1){#3}\hss}}
\def\Pfad(#1,#2),#3\endPfad{\unskip\leavevmode
  \xcoord#1 \ycoord#2 \thicklines\ZeichnePfad#3\endPfad\thinlines}
\def\ZeichnePfad#1{\ifx#1\endPfad\let\next\relax
  \else\let\next\ZeichnePfad
    \ifnum#1=1
      \raise\ycoord \Einheit\hbox to0pt{\hskip\xcoord \Einheit
         \vrule height\Pfadd@cke width1 \Einheit depth\Pfadd@cke\hss}%
      \advance\xcoord by 1
    \else\ifnum#1=2
      \raise\ycoord \Einheit\hbox to0pt{\hskip\xcoord \Einheit
        \hbox{\hskip-\PfadD@cke\vrule height1 \Einheit width\PfadD@cke depth0pt}\hss}%
      \advance\ycoord by 1
    \else\ifnum#1=3
      \raise\ycoord \Einheit\hbox to0pt{\hskip\xcoord \Einheit
         \Line@(1,1){1}\hss}
      \advance\xcoord by 1
      \advance\ycoord by 1
    \else\ifnum#1=4
      \raise\ycoord \Einheit\hbox to0pt{\hskip\xcoord \Einheit
         \Line@(1,-1){1}\hss}
      \advance\xcoord by 1
      \advance\ycoord by -1
    \else\ifnum#1=5
      \advance\xcoord by -1
      \raise\ycoord \Einheit\hbox to0pt{\hskip\xcoord \Einheit
         \vrule height\Pfadd@cke width1 \Einheit depth\Pfadd@cke\hss}%
    \else\ifnum#1=6
      \advance\ycoord by -1
      \raise\ycoord \Einheit\hbox to0pt{\hskip\xcoord \Einheit
        \hbox{\hskip-\PfadD@cke\vrule height1 \Einheit width\PfadD@cke depth0pt}\hss}%
    \else\ifnum#1=7
      \advance\xcoord by -1
      \advance\ycoord by -1
      \raise\ycoord \Einheit\hbox to0pt{\hskip\xcoord \Einheit
         \Line@(1,1){1}\hss}
    \else\ifnum#1=8
      \advance\xcoord by -1
      \advance\ycoord by +1
      \raise\ycoord \Einheit\hbox to0pt{\hskip\xcoord \Einheit
         \Line@(1,-1){1}\hss}
    \fi\fi\fi\fi
    \fi\fi\fi\fi
  \fi\next}
\def\hSSchritt{\leavevmode\raise-.4pt\hbox to0pt{\hss.\hss}\hskip.2\Einheit
  \raise-.4pt\hbox to0pt{\hss.\hss}\hskip.2\Einheit
  \raise-.4pt\hbox to0pt{\hss.\hss}\hskip.2\Einheit
  \raise-.4pt\hbox to0pt{\hss.\hss}\hskip.2\Einheit
  \raise-.4pt\hbox to0pt{\hss.\hss}\hskip.2\Einheit}
\def\vSSchritt{\vbox{\baselineskip.2\Einheit\lineskiplimit0pt
\hbox{.}\hbox{.}\hbox{.}\hbox{.}\hbox{.}}}
\def\DSSchritt{\leavevmode\raise-.4pt\hbox to0pt{%
  \hbox to0pt{\hss.\hss}\hskip.2\Einheit
  \raise.2\Einheit\hbox to0pt{\hss.\hss}\hskip.2\Einheit
  \raise.4\Einheit\hbox to0pt{\hss.\hss}\hskip.2\Einheit
  \raise.6\Einheit\hbox to0pt{\hss.\hss}\hskip.2\Einheit
  \raise.8\Einheit\hbox to0pt{\hss.\hss}\hss}}
\def\dSSchritt{\leavevmode\raise-.4pt\hbox to0pt{%
  \hbox to0pt{\hss.\hss}\hskip.2\Einheit
  \raise-.2\Einheit\hbox to0pt{\hss.\hss}\hskip.2\Einheit
  \raise-.4\Einheit\hbox to0pt{\hss.\hss}\hskip.2\Einheit
  \raise-.6\Einheit\hbox to0pt{\hss.\hss}\hskip.2\Einheit
  \raise-.8\Einheit\hbox to0pt{\hss.\hss}\hss}}
\def\SPfad(#1,#2),#3\endSPfad{\unskip\leavevmode
  \xcoord#1 \ycoord#2 \ZeichneSPfad#3\endSPfad}
\def\ZeichneSPfad#1{\ifx#1\endSPfad\let\next\relax
  \else\let\next\ZeichneSPfad
    \ifnum#1=1
      \raise\ycoord \Einheit\hbox to0pt{\hskip\xcoord \Einheit
         \hSSchritt\hss}%
      \advance\xcoord by 1
    \else\ifnum#1=2
      \raise\ycoord \Einheit\hbox to0pt{\hskip\xcoord \Einheit
        \hbox{\hskip-2pt \vSSchritt}\hss}%
      \advance\ycoord by 1
    \else\ifnum#1=3
      \raise\ycoord \Einheit\hbox to0pt{\hskip\xcoord \Einheit
         \DSSchritt\hss}
      \advance\xcoord by 1
      \advance\ycoord by 1
    \else\ifnum#1=4
      \raise\ycoord \Einheit\hbox to0pt{\hskip\xcoord \Einheit
         \dSSchritt\hss}
      \advance\xcoord by 1
      \advance\ycoord by -1
    \else\ifnum#1=5
      \advance\xcoord by -1
      \raise\ycoord \Einheit\hbox to0pt{\hskip\xcoord \Einheit
         \hSSchritt\hss}%
    \else\ifnum#1=6
      \advance\ycoord by -1
      \raise\ycoord \Einheit\hbox to0pt{\hskip\xcoord \Einheit
        \hbox{\hskip-2pt \vSSchritt}\hss}%
    \else\ifnum#1=7
      \advance\xcoord by -1
      \advance\ycoord by -1
      \raise\ycoord \Einheit\hbox to0pt{\hskip\xcoord \Einheit
         \DSSchritt\hss}
    \else\ifnum#1=8
      \advance\xcoord by -1
      \advance\ycoord by 1
      \raise\ycoord \Einheit\hbox to0pt{\hskip\xcoord \Einheit
         \dSSchritt\hss}
    \fi\fi\fi\fi
    \fi\fi\fi\fi
  \fi\next}
\def\Koordinatenachsen(#1,#2){\unskip
 \hbox to0pt{\hskip-.5pt\vrule height#2 \Einheit width.5pt depth1 \Einheit}%
 \hbox to0pt{\hskip-1 \Einheit \xcoord#1 \advance\xcoord by1
    \vrule height0.25pt width\xcoord \Einheit depth0.25pt\hss}}
\def\Koordinatenachsen(#1,#2)(#3,#4){\unskip
 \hbox to0pt{\hskip-.5pt \ycoord-#4 \advance\ycoord by1
    \vrule height#2 \Einheit width.5pt depth\ycoord \Einheit}%
 \hbox to0pt{\hskip-1 \Einheit \hskip#3\Einheit 
    \xcoord#1 \advance\xcoord by1 \advance\xcoord by-#3 
    \vrule height0.25pt width\xcoord \Einheit depth0.25pt\hss}}
\def\Gitter(#1,#2){\unskip \xcoord0 \ycoord0 \leavevmode
  \LOOP\ifnum\ycoord<#2
    \loop\ifnum\xcoord<#1
      \raise\ycoord \Einheit\hbox to0pt{\hskip\xcoord \Einheit\Punkt\hss}%
      \advance\xcoord by1
    \repeat
    \xcoord0
    \advance\ycoord by1
  \REPEAT}
\def\Gitter(#1,#2)(#3,#4){\unskip \xcoord#3 \ycoord#4 \leavevmode
  \LOOP\ifnum\ycoord<#2
    \loop\ifnum\xcoord<#1
      \raise\ycoord \Einheit\hbox to0pt{\hskip\xcoord \Einheit\Punkt\hss}%
      \advance\xcoord by1
    \repeat
    \xcoord#3
    \advance\ycoord by1
  \REPEAT}
\def\Label#1#2(#3,#4){\unskip \xdim#3 \Einheit \ydim#4 \Einheit
  \def\lo{\advance\xdim by-.5 \Einheit \advance\ydim by.5 \Einheit}%
  \def\llo{\advance\xdim by-.25cm \advance\ydim by.5 \Einheit}%
  \def\loo{\advance\xdim by-.5 \Einheit \advance\ydim by.25cm}%
  \def\o{\advance\ydim by.25cm}%
  \def\ro{\advance\xdim by.5 \Einheit \advance\ydim by.5 \Einheit}%
  \def\rro{\advance\xdim by.25cm \advance\ydim by.5 \Einheit}%
  \def\roo{\advance\xdim by.5 \Einheit \advance\ydim by.25cm}%
  \def\l{\advance\xdim by-.30cm}%
  \def\r{\advance\xdim by.30cm}%
  \def\lu{\advance\xdim by-.5 \Einheit \advance\ydim by-.6 \Einheit}%
  \def\llu{\advance\xdim by-.25cm \advance\ydim by-.6 \Einheit}%
  \def\luu{\advance\xdim by-.5 \Einheit \advance\ydim by-.30cm}%
  \def\u{\advance\ydim by-.30cm}%
  \def\ru{\advance\xdim by.5 \Einheit \advance\ydim by-.6 \Einheit}%
  \def\rru{\advance\xdim by.25cm \advance\ydim by-.6 \Einheit}%
  \def\ruu{\advance\xdim by.5 \Einheit \advance\ydim by-.30cm}%
  #1\raise\ydim\hbox to0pt{\hskip\xdim
     \vbox to0pt{\vss\hbox to0pt{\hss$#2$\hss}\vss}\hss}%
}
\catcode`\@=13

\catcode`\@=11
\def\logo@{}
\footline={\ifnum\pageno>1 \hfil\folio\hfil\else\hfil\fi}
\topmatter
\title Enumeration of lozenge tilings of hexagons with cut off corners 
\endtitle
\author Mihai Ciucu {\rm and} Christian Krattenthaler\endauthor
\subjclass
Primary 05A15, 05A17, 05B45. 
Secondary 11P81, 52C20
\endsubjclass
\keywords
Plane partitions, symmetry classes, determinant evaluations, lozenge tilings, 
non-intersecting lattice paths, tiling enumeration, perfect matchings
\endkeywords

\thanks The first author was partially supported by NSF grant DMS-9802390.
\endthanks
\thanks The second author was partially supported by the Austrian
Science Foundation FWF, grant P13190-MAT.
\endthanks
\affil
  Georgia Institute of Technology\\
  School of Mathematics\\
  Atlanta, GA 30332, USA\\
\\
and\\
\\
Institut f\"ur Mathematik der Universit\"at Wien\\
Strudlhofgasse 4, A-1090 Wien, Austria
\endaffil
\date February, 2000\enddate
\abstract
Motivated by the enumeration of a class of plane partitions studied by Proctor and
by considerations about symmetry classes of plane partitions, we consider
the problem of enumerating lozenge tilings of a hexagon with ``maximal staircases'' removed from 
some of its vertices. The case of one vertex corresponds to Proctor's
problem. For two vertices there are several cases to consider, and most of them lead
to nice enumeration formulas. 
For three or more vertices there do not seem to exist
nice product formulas in general, but in one special situation a lot of factorization
occurs, and we pose the problem of finding a formula for the number of tilings in
this case.
\endabstract
\endtopmatter
\document

\def\mysec#1{\bigskip\centerline{\bf #1}\message{ * }\nopagebreak\par}

\def\myref#1{\item"{[{\bf #1}]}"} 
 
\def\pf{{\it Proof.\ }} 

\def\cite#1{\relaxnext@
  \def\nextiii@##1,##2\end@{[{\bf##1},\,##2]}%
  \in@,{#1}\ifin@\def\next{\nextiii@#1\end@}\else
  \def\next{[{\bf#1}]}\fi\next}
\def\proclaimheadfont@{\smc}

\def\pf{{\it Proof.\ }}

\define\Z{{\bold Z}}

\define\twoline#1#2{\line{\hfill{\smc #1}\hfill{\smc #2}\hfill}}
\define\twolinetwo#1#2{\line{{\smc #1}\hfill{\smc #2}}}
\define\twolinethree#1#2{\line{\phantom{poco}{\smc #1}\hfill{\smc #2}\phantom{poco}}}

\def\mypic#1{\epsffile{#1}}
%\def\mypic#1{\epsffile{#1}}

%references

\def\Stenlp{25}
\def\Stapp{24}
\def\Robb{23}
\def\Proc88{22}
\def\PeWZAA{21}
\def\MacM{20}
\def\LindAA{19}
\def\Ku{18}
\def\KratBN{17}
\def\KratBD{16}
\def\Kast{15}
\def\GrKPAA{14}
\def\GospAB{13}
\def\GV{12}
\def\FuKrAC{11}
\def\Feller{10}
\def\EisTAA{9}
\def\DT{8}
\def\CiKrAD{7}
\def\CKpp2{6}
\def\CiKrAA{5}
\def\Cipp1{4}
\def\CiucAI{3}
\def\Ci1{2}
\def\Amdeb{1}

\def\({\left(}
\def\){\right)}
\def\[{\left[}
\def\]{\right]}

\mysec{1. Introduction and statement of results}

\medskip
The study of lozenge tilings is warranted by the many useful insights they bring in
the subject of plane partitions. Some important instances of these are presented in
\cite{\Ku} and \cite{\CKpp2}. In this paper we present some more such connections.

%\medskip
A plane partition is a rectangular array of nonnegative integers with the property 
that all rows
and columns are weakly decreasing. A plane partition contained in an $a\times b$ 
rectangle and with entries at most $c$ can be identified with its three dimensional
diagram --- a stack of unit cubes contained in an $a\times b\times c$ box --- ,
which in turn can be regarded as a lozenge tiling of a hexagon $H(a,b,c)$ with 
side lengths $a$, $b$,
$c$, $a$, $b$, $c$ (in cyclic order) and angles of $120^\circ$ (see Figure~1.1
and \cite{\DT} or \cite{\Robb}; a lozenge tiling of
a region on the triangular lattice is a tiling by unit rhombi with angles of
$60^\circ$ and $120^\circ$, referred to as lozenges). 
This simple bijection is the crucial link between the theory of lozenge tilings and
that of plane partitions. For example, the number of tilings of $H(a,b,c)$ follows to
be equal to the number of plane partitions that fit in an $a\times b\times c$ box,
which is, 
by a result due to MacMahon \cite{\MacM},
$\prod_{i=1}^a\prod_{j=1}^b\prod_{k=1}^c 
(i+j+k-1)/(i+j+k-2)$.

%\medskip
As a variation of this, Proctor \cite{\Proc88} considered the problem of enumerating
those plane partitions $\pi$ contained in an $a\times b\times c$ box for which the
projection of $\pi$ on one of the coordinate planes, say on $Oxy$, fits in the ``maximal
staircase'' $\lambda=(b,b-1,\dotsc,b-a+1)$ 
(when $\lambda$ is viewed as the corresponding Ferrers diagram) 
contained in the $a\times b$ basis 
rectangle (we are assuming here, without loss of generality, that $a\leq b;$ see
Figure~1.1 for an example of such a plane partition with $a=5$, $b=8$, $c=3$). 
Proctor  \cite{\Proc88} found that this number is given by the simple product
$$\prod_{i=1}^a\left[\prod_{j=1}^{b-a+1}\frac{c+i+j-1}{i+j-1}
\prod_{j=b-a+2}^{b-a+i}\frac{2c+i+j-1}{i+j-1}\right].\tag1.1$$
By the above bijection, it is easily seen that Proctor's
problem is equivalent to counting the lozenge tilings of $H(a,b,c)$ with a maximal
``staircase of lozenges'' removed from a corner at which edges of lengths $a$ and $b$
meet (simply view Figure~1.1 as being two dimensional%
; there, a maximal staircase of lozenges was removed from the
  southeastern corner). 

%\medskip
What if we require that the projection of $\pi$ on {\it two} of the coordinate
planes be contained in the corresponding staircases? The above bijection shows that
the question is equivalent to counting the number of tilings of $H(a,b,c)$ with 
{\it two} maximal staircases removed, from vertices that are non-adjacent and
non-opposite (this is illustrated in Figure~1.2; 
there, maximal staircases were removed from the
  southeastern and the western corner; the plane partition is the same as
in Figure~1.1). There are six cases to consider, corresponding to the six relative
orderings of $a$, $b$ and $c$. These are shown in Figures~1.3(a)--(f). 
(At this point, the special marks in form of ellipses should be
ignored.)
Mirror reflection pairs up these six cases in three pairs --- the rows of Figure~1.3. 

\topinsert
\twoline{\mypic{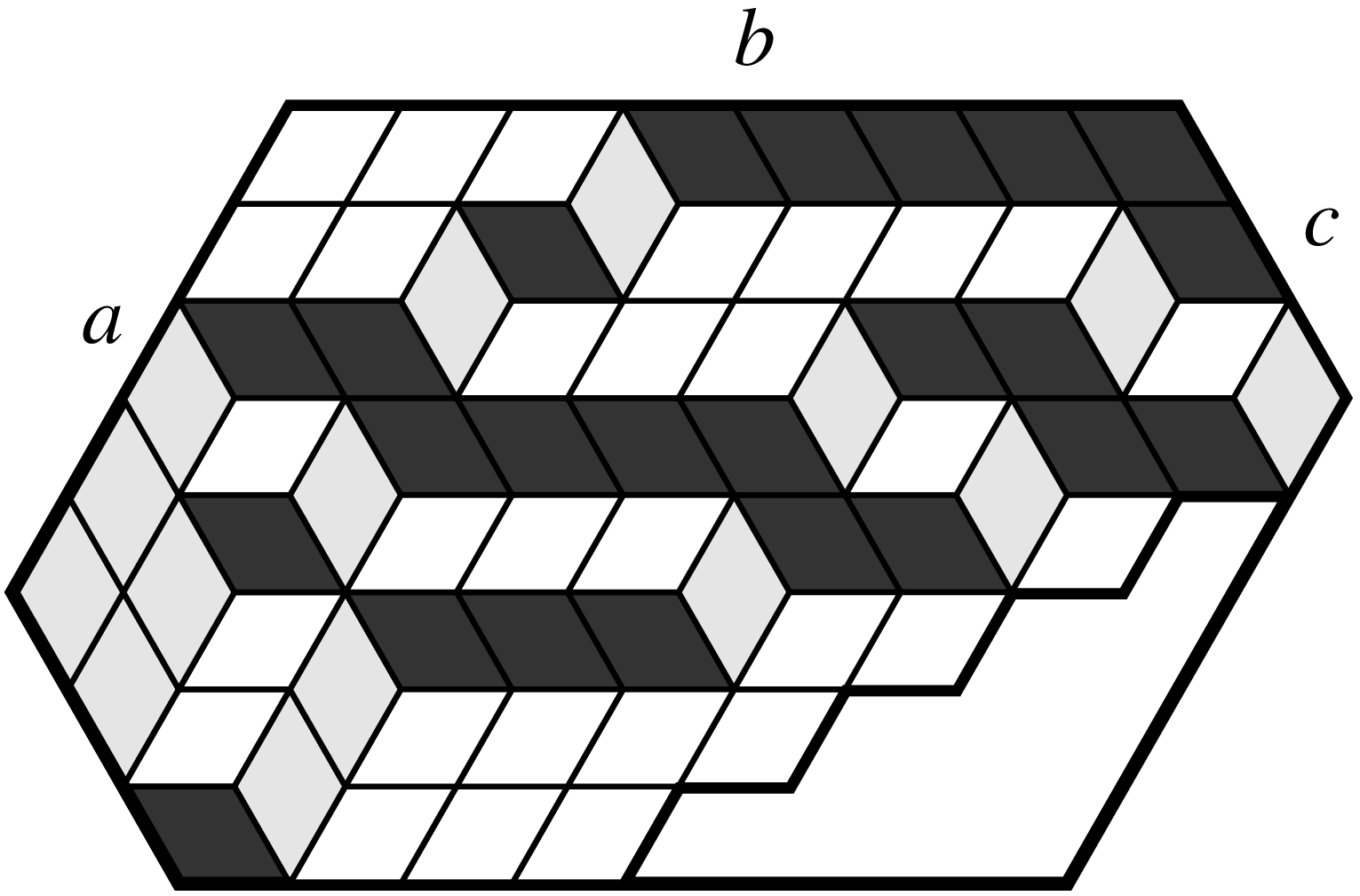}}{\mypic{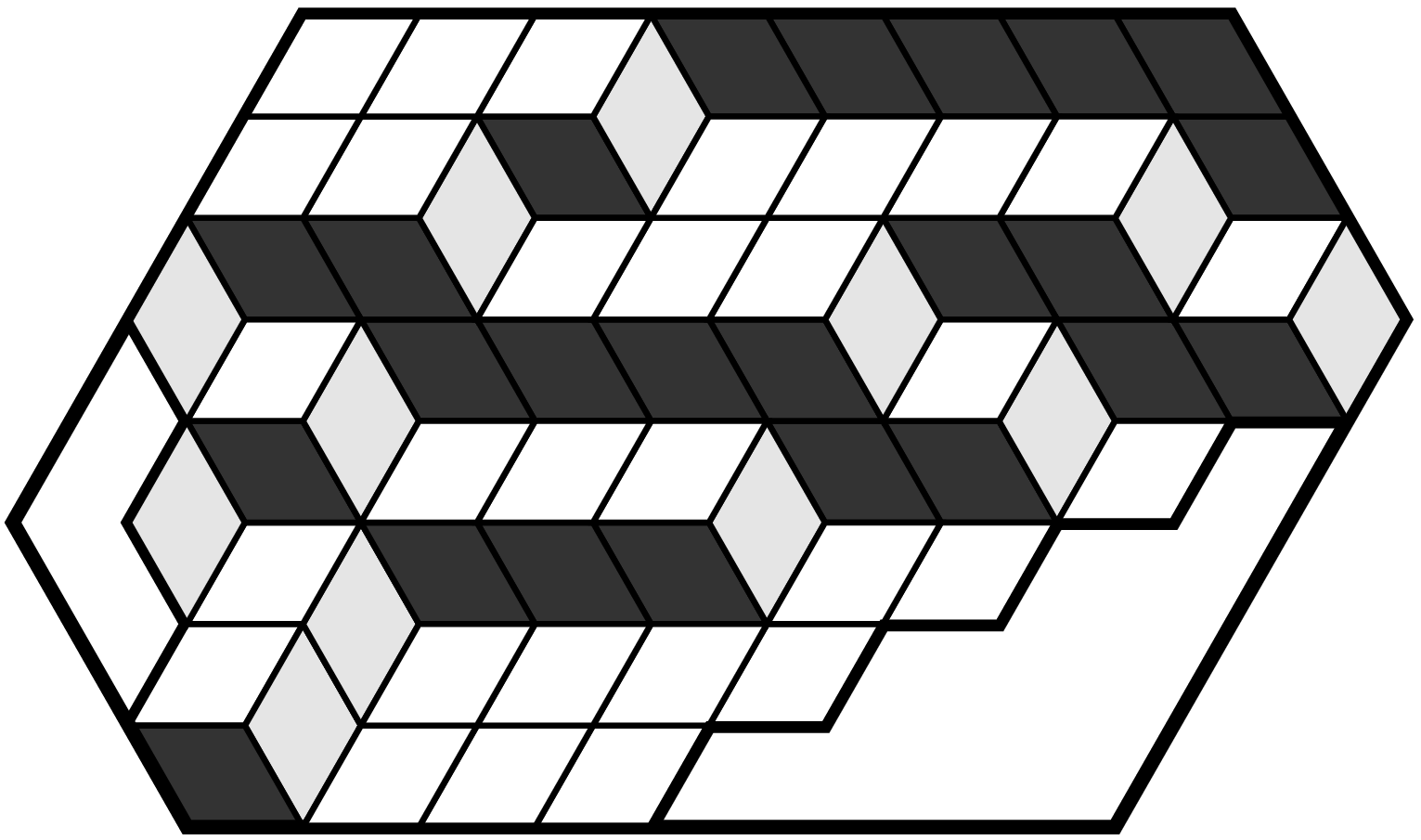}}
\twoline{Figure~1.1. {\rm \ \ \ \ \  }}{\ \ \ \ \ \ \ \ \ \ Figure~1.2. {\rm }}
\endinsert

\topinsert
\twoline{\mypic{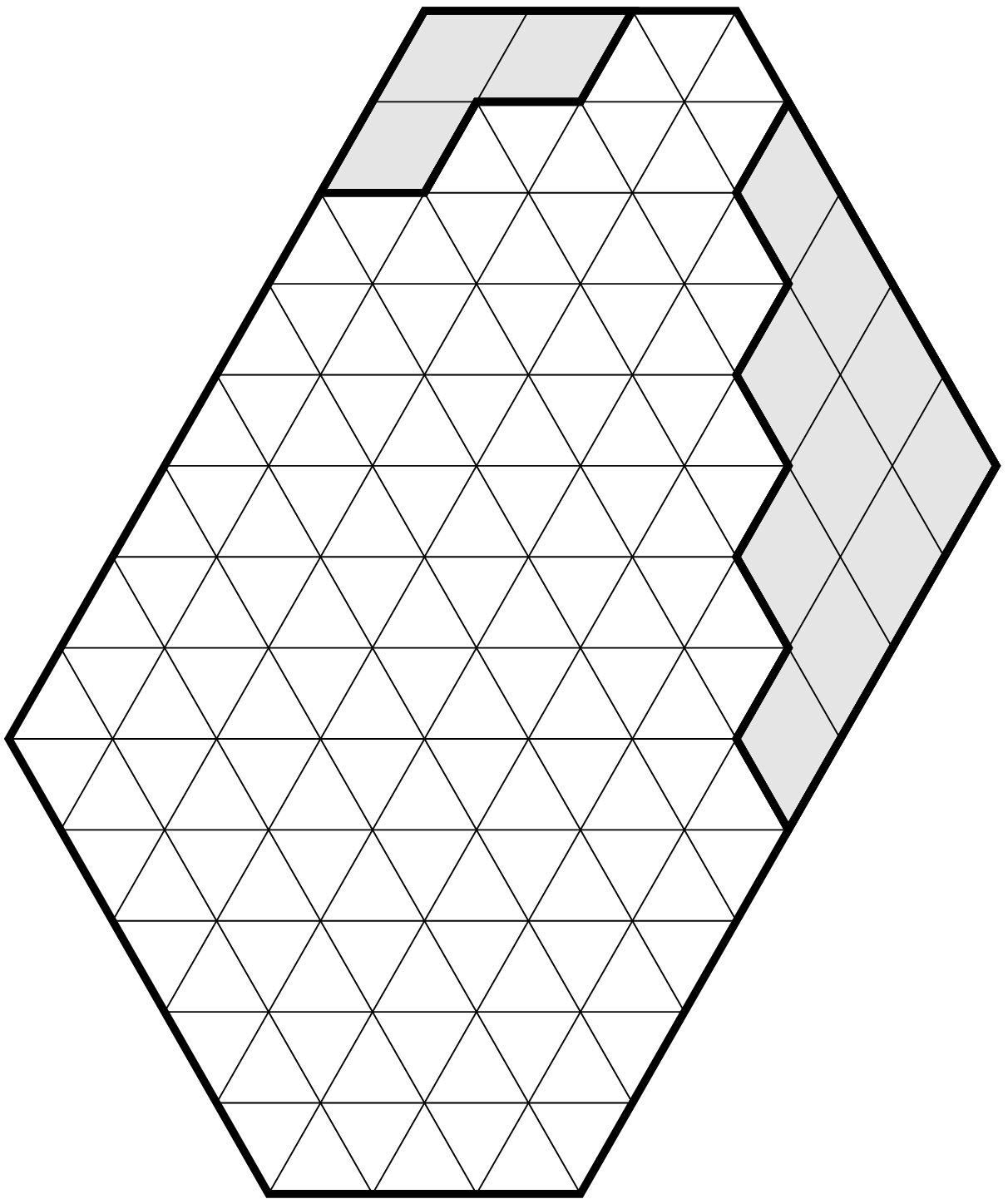}}{\mypic{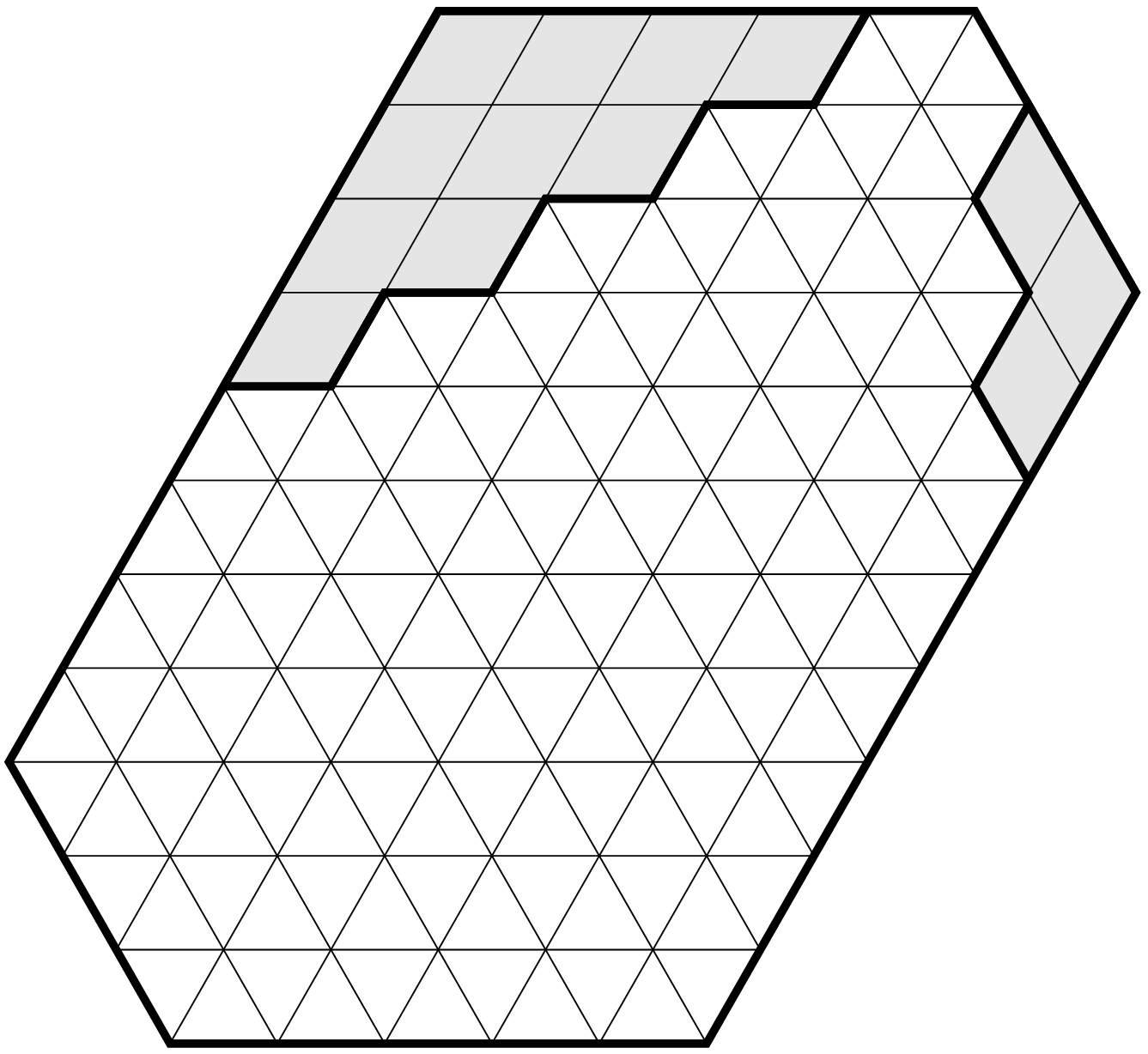}}
\twoline{Figure~1.3{\rm (a). $H_{\text d}(8,3,5)$.}}
{Figure~1.3{\rm (b). $H_{\text d}(8,5,3)$.}}
\endinsert

\topinsert
\twoline{\mypic{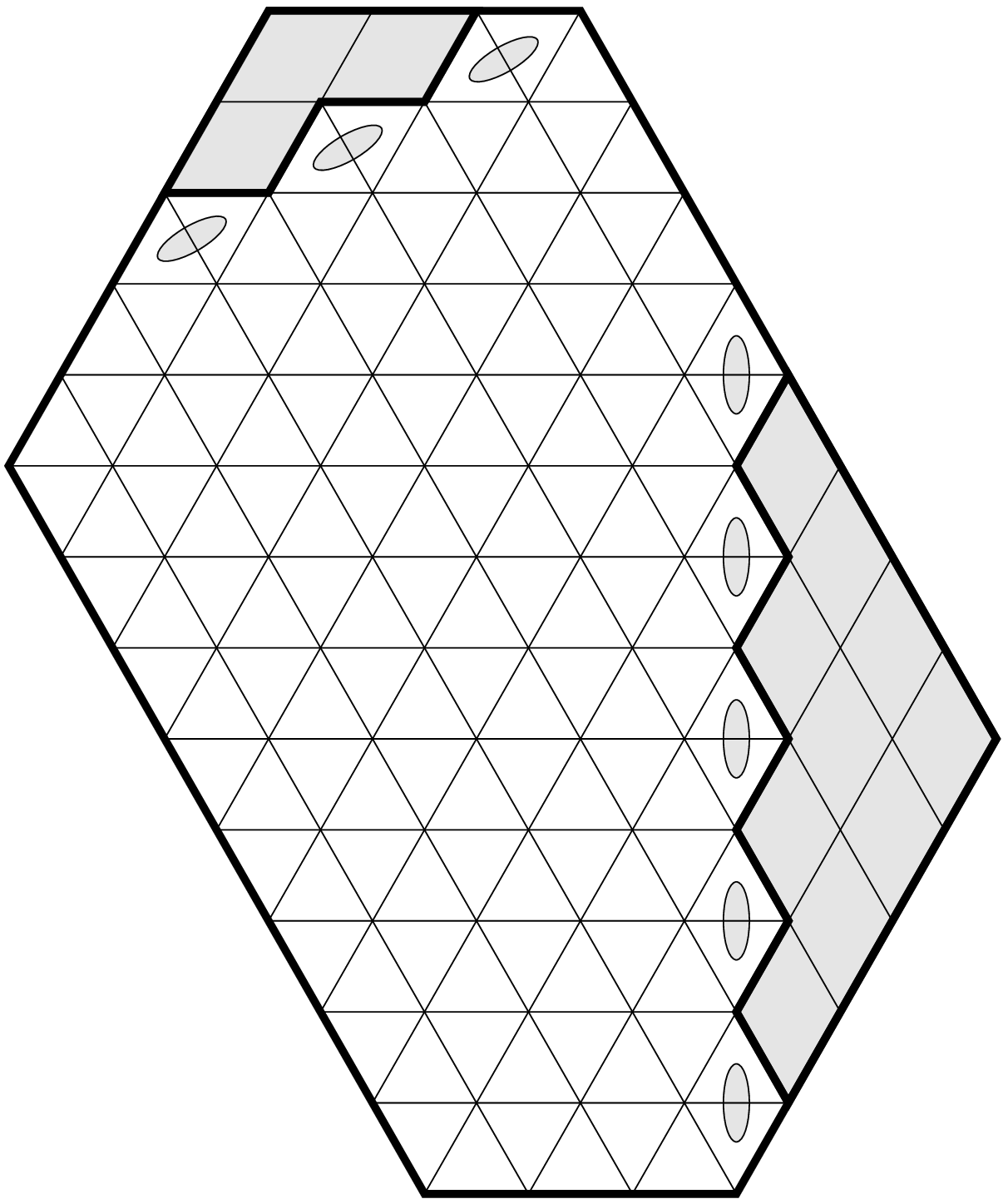}}{\mypic{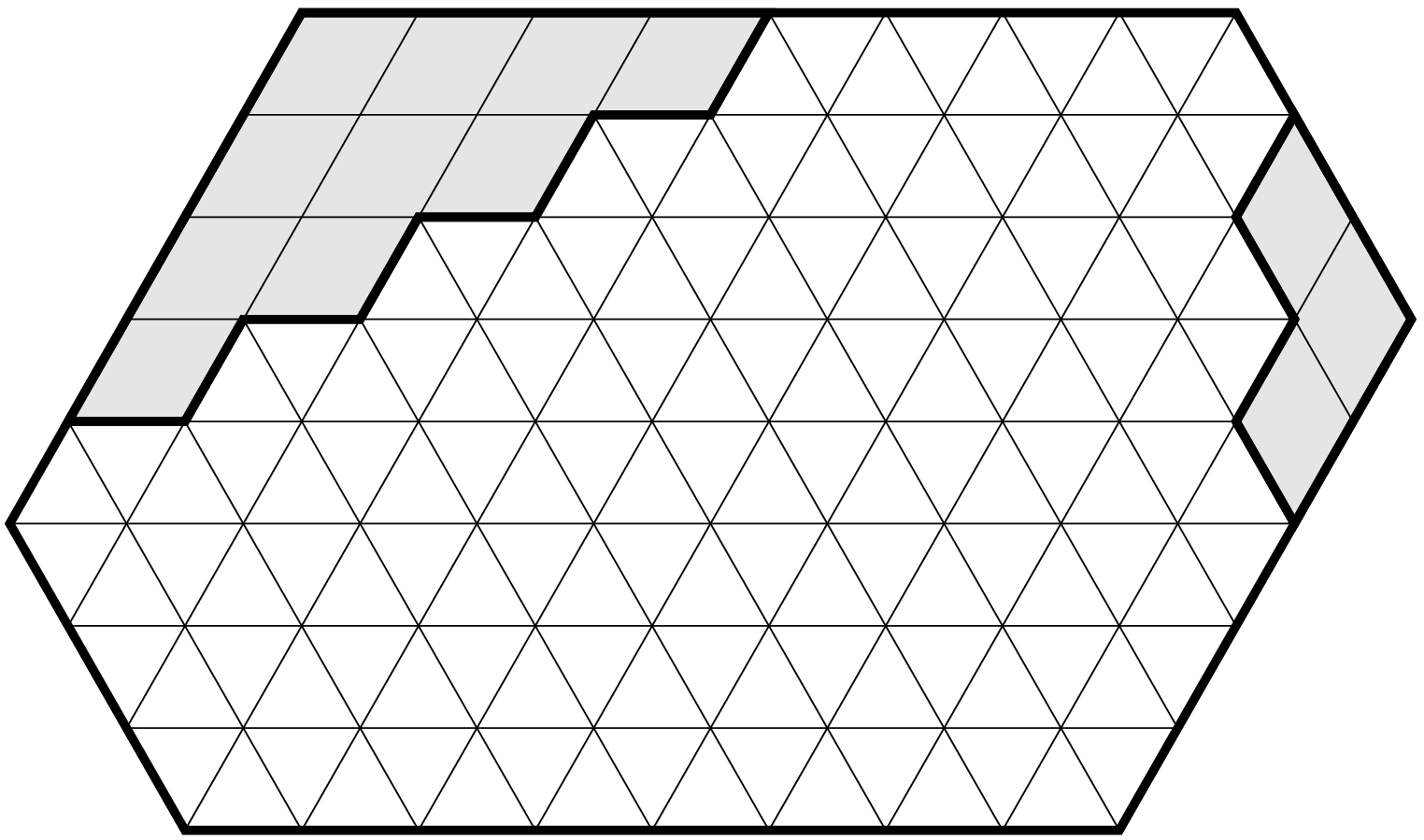}}
\twoline{Figure~1.3{\rm (c). $H_{\text d}(5,3,8)$.}}
{Figure~1.3{\rm (d). $H_{\text d}(5,8,3)$.}}
\endinsert

\topinsert
\twoline{\mypic{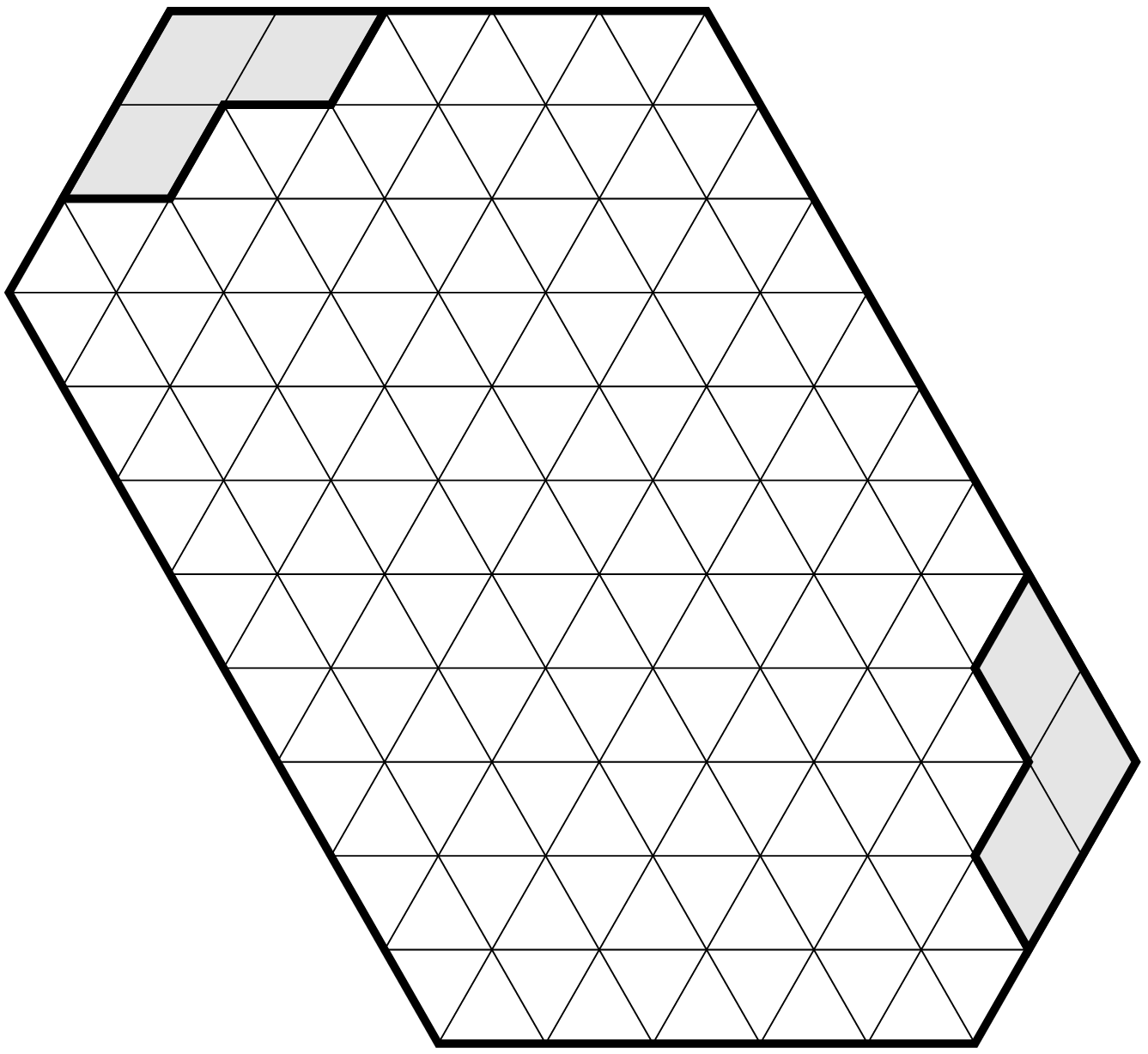}}{\mypic{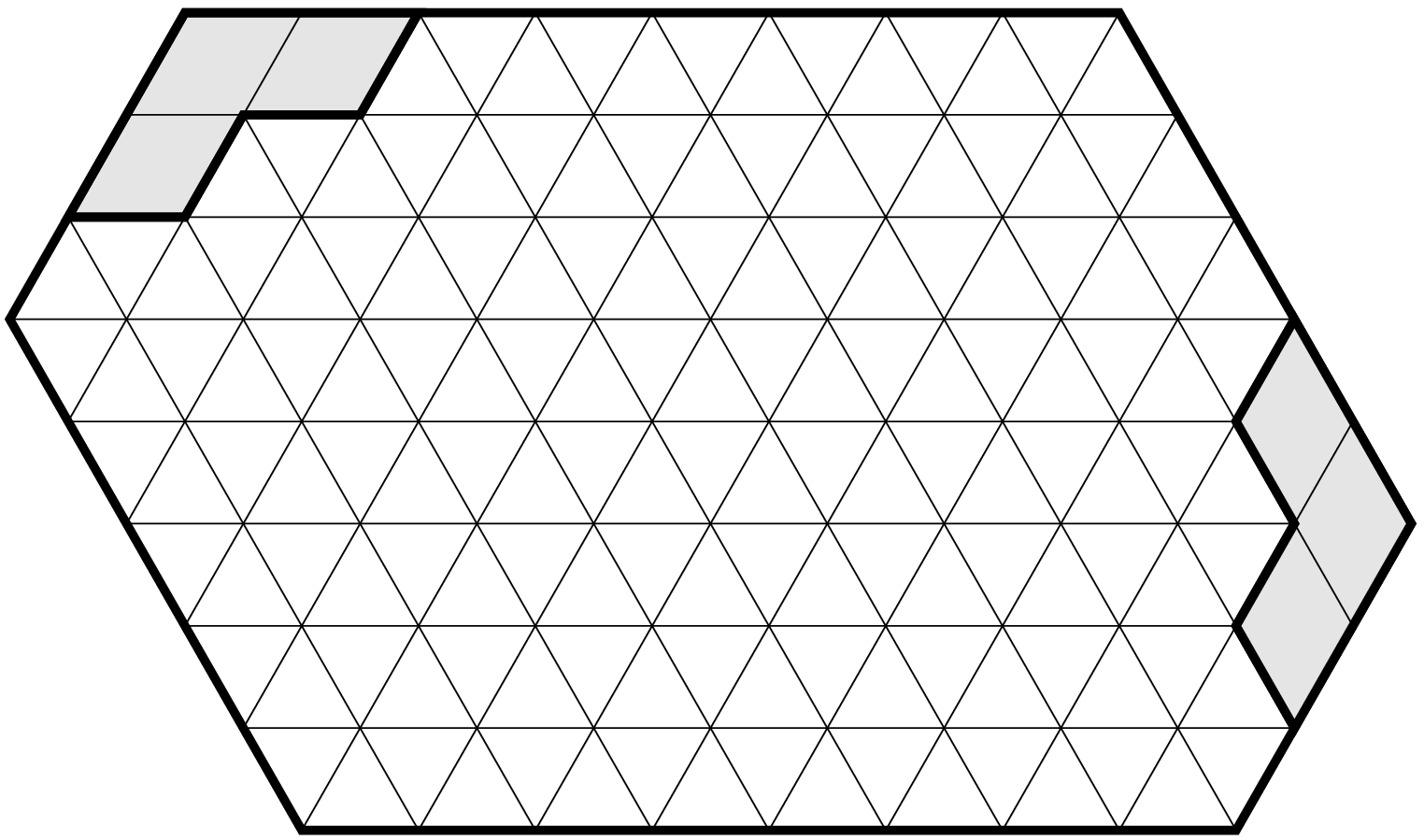}}
\twoline{Figure~1.3{\rm (e). $H_{\text d}(3,5,8)$.}}
{Figure~1.3{\rm (f). $H_{\text d}(3,8,5)$.}}
\endinsert

%\medskip
We draw all the hexagons $H(a,b,c)$ and the regions obtained from them by
cutting corners so that the horizontal edges have length $b$, and the other two pairs
of parallel edges, as we move counterclockwise, have lengths $a$ and $c$.  
Let $H_{\text d}(a,b,c)$ be the region obtained from $H(a,b,c)$ by removing
maximal staircases from the northwestern and eastern corners (the subscript stands for 
the ``diagonal'' position of the cut off corners).  

%Define M, M_*, M^*, M_*^*.

%\medskip
For a region $R$ on the triangular lattice, denote by $L(R)$ the number of
its lozenge tilings. In the special case when $R$ is obtained from a hexagon
$H(a,b,c)$ by removing staircases from two of its corners, we define three more
{\it weighted\/} tiling enumerators for $R$ as follows. Consider the tile positions that fit
in the indentations of the zig-zag cut that removed a staircase of lozenges (the
possible such positions are marked 
by ellipses in Figure~1.3(c)). By
weighting these tile positions by 1/2, one creates a new, weighted count of the
tilings of $R$: each tiling $T$ gets weight $1/2^k$, where $k$ is the number of 
lozenges of $T$ occupying positions weighted by 1/2, and the sum of weights of all
tilings $T$ of $R$ gives the weighted tiling enumeration
\footnote{The motivation to consider such weightings comes from the
fact that weightings of that kind arise whenever the Factorization
Theorem from \cite{\Ci1} is applied to a (symmetric) region on the
triangular lattice. 
See \cite{\Ci1, Sec.~6}\cite{\CiucAI}\cite{\Cipp1}\cite{\CiKrAA}\cite{\CKpp2}%
\cite{\EisTAA}\cite{\FuKrAC} and the proof of Theorem~1.4 in Section~2
for examples.}.
Clearly, one can choose to
weight by 1/2 {\it only} the tile positions along the cut that removed the 
northwestern corner of 
$H(a,b,c)$, or, furthermore, to weight by 1/2 {\it only} the tile
positions along the cut that removed the  
eastern corner. These three
possibilities define our three weighted enumerators. We denote them by
$L_*^*$, $L^*$, and $L_*$, where a superscript (respectively,
subscript) indicates weighting along the cut  
from the northwestern (respectively, eastern) corner.

If $b<c<a$ (see Figure~1.3(a)) --- or, by mirror reflection, $c<b<a$ 
(see Figure~1.3(b)), --- neither $L(H_{\text d}(a,b,c))$, nor the
weighted enumerators $L^*(H_{\text d}(a,b,c))$, $L_*(H_{\text
d}(a,b,c))$, and $L^*_*(H_{\text d}(a,b,c))$
seem to be given by simple product formulas. The other
two cases lead to the following results.

\proclaim{Theorem 1.1} If $b\leq a\leq c$ $($see Figure~{\rm1.3(c)} for an
example$)$, we have
$$L(H_{\text d}(a,b,c))=(-1)^{b(b+1)/2} P_b(a-2b-1,b+c+1),\tag1.2$$
where $P_n(x,y)$ denotes the product on the right hand side of\/ 
{\rm(1.8)}.
\endproclaim

\flushpar
{\bf Note.} 
All the factors in $P_b(a-2b-1,b+c+1)$ are positive
except for the factors in the shifted factorial $(a-3b-c+2j-1)_j$, which are all negative
since for the largest factor in this product we have
$$a-3b-c+3j-2\leq a-c-2\leq -2,$$
as $a\leq c$ in the case addressed by Theorem~1.1. Therefore, for $a$, $b$ and $c$ as
in Theorem~1.1, the sign of $P_b(a-2b-1,b+c+1)$ is $(-1)^{b(b+1)/2}$, which checks that
the right hand side of (1.2) is non-negative.

Still keeping the relative order $b\le a\le c$ of the parameters,
the weighted enumerators $L^*(H_{\text d}(a,b,c))$
and $L_*(H_{\text d}(a,b,c))$ 
do not seem to be given by simple product formulas. But there is one
for $L^*_*(H_{\text d}(a,b,c))$.

\proclaim{Theorem 1.2} If $b\leq a\leq c$ $($see Figure~{\rm1.3(c)} for an
example$)$, we have
$$\multline
L^*_*(H_{\text d}(a,b,c))\\
=2^{-a-b}
\prod _{j=1} ^{b}\frac {(j-1)!\, (a+c-b+2j-1)! \,(c-a+3j-1)_{b-j}\,
(a+2c+3j-1)_{b-j+1}} {(b+c+j-1)! \,(a-b+2j-1)!},
\endmultline\tag1.3$$
where the shifted factorial
$(\alpha)_k$ is defined by $(\alpha)_k:=\alpha(\alpha+1)\cdots(\alpha+k-1)$, 
$k\ge1$, and $(\alpha)_0:=1$.
\endproclaim

\smallpagebreak
For $a\leq b\leq c$, plain enumeration of the tilings of $H_{\text d}(a,b,c)$
does not seem to be given by a simple product formula. There is, however, a simple
formula for $L^*(H_{\text d}(a,b,c))$.
As we are going to show in Section~2, the following result follows
easily from a determinant evaluation of the second author 
\cite{\KratBD, (5.3)}. 

\proclaim{Proposition 1.3} If $a\leq b\leq c$ $($see
Figure~{\rm1.3(e)} for an
example$)$, 
we have
$$\multline
L^*(H_{\text d}(a,b,c))=\frac{1}{2^a}\prod_{j=1}^{a}
\frac{(j-1)!\,(b+2c-3a+3j)\,(b+c-2a+j)!}
{(c-a+2j-1)!}\\
\cdot\frac {(b+2c-3a+2j+1)_{j-1}\,(2b+c-3a+2j+1)_{j-1}}
{(b-a+2j-1)!}.
\endmultline$$

%L_*(H_{\text d}(a,b,c))=... is equivalent to above by swapping b and c

\endproclaim

\flushpar
{\bf Note.} The special cases when $c=a-1$ or when $b=a-1$ form the
subject of \cite{\CKpp2, Pro\-pos\-i\-tion~2.2}.

\smallpagebreak
The weighted enumerator $L^*$ clearly makes sense also for the region $H_1(a,b,c)$
obtained from $H(a,b,c)$ by cutting off just the northwestern corner. We have the
following counterpart of Proctor's formula (1.1).

\proclaim{Theorem 1.4} For $a\leq b$, we have
$$L^*(H_1(a,b,c))=\frac{1}{2^a}\prod_{i=1}^a\left[\prod_{j=1}^{b-a}\frac{c+i+j-1}{i+j-1}
\prod_{j=b-a+i}^{b}\frac{2c+i+j-1}{i+j-1}\right].\tag1.4$$
\endproclaim

What if we require that the projection of $\pi$ on {\it all three} coordinate planes be
contained in the corresponding staircases? An example illustrating this, using the
same plane partition as in Figures~1.1 and 1.2, is shown in Figure~1.4.
Clearly, by the bijection between plane partitions and lozenge tilings, plane
partitions satisfying these conditions are identified with
tilings of the region $H_3(a,b,c)$ obtained from $H(a,b,c)$ by removing maximal
staircases from three alternating vertices.
No matter what the relative ordering of the side lengths $a$, $b$ and $c$ is, 
there are always two staircases that
``interfere'' --- i.e., there is no portion of an edge of the hexagon $H(a,b,c)$
separating them. 
Data suggests that there are no simple product formulas in
general. 
(Note that, among the cases with two removed staircases, the relative
orders of $a,b,c$ not covered by Theorems~1.1, 1.2
and Proposition~1.3 are precisely those in which the staircases interfere.) 

However, the special case $a=b=c$ presents significant factorization. Indeed, letting 
$T_a$ stand for the triangular region $H_3(a,a,a)$ (the case $a=5$ is shown in
Figure~1.5), the number of tilings of $T_a$ factors as follows for $a\leq7$:
$$\align
L(T_1)&=2\\ 
L(T_2)&=3^2\\ 
L(T_3)&=2^3\cdot13\\ 
L(T_4)&=2^2\cdot5^2\cdot31\\ 
L(T_5)&=2\cdot3^2\cdot19^2\cdot37\\ 
L(T_6)&=2\cdot7^3\cdot13\cdot43\cdot127\\ 
L(T_7)&=2^7\cdot3^5\cdot5^3\cdot7\cdot13\cdot73.
\endalign$$ 
The amount of factorization is
remarkable (also for larger $a$; we have computed and factored $L(T_a)$ up to
$a=30$) and
comparable, say, to that of the
numbers enumerating domino tilings of squares (see \cite{\Kast}). Based on this, we
pose the following problem.

\topinsert
\twoline{\mypic{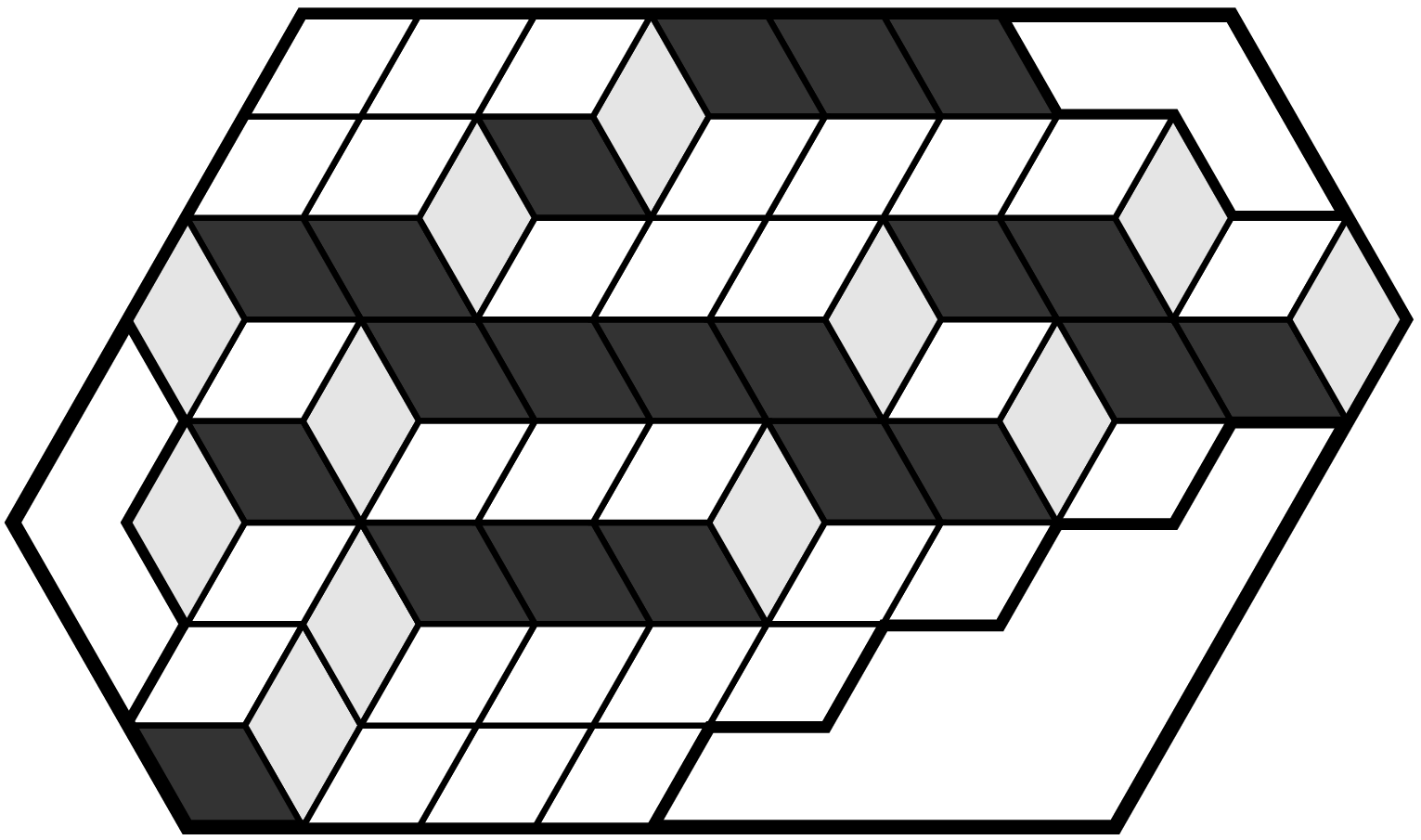}}{\mypic{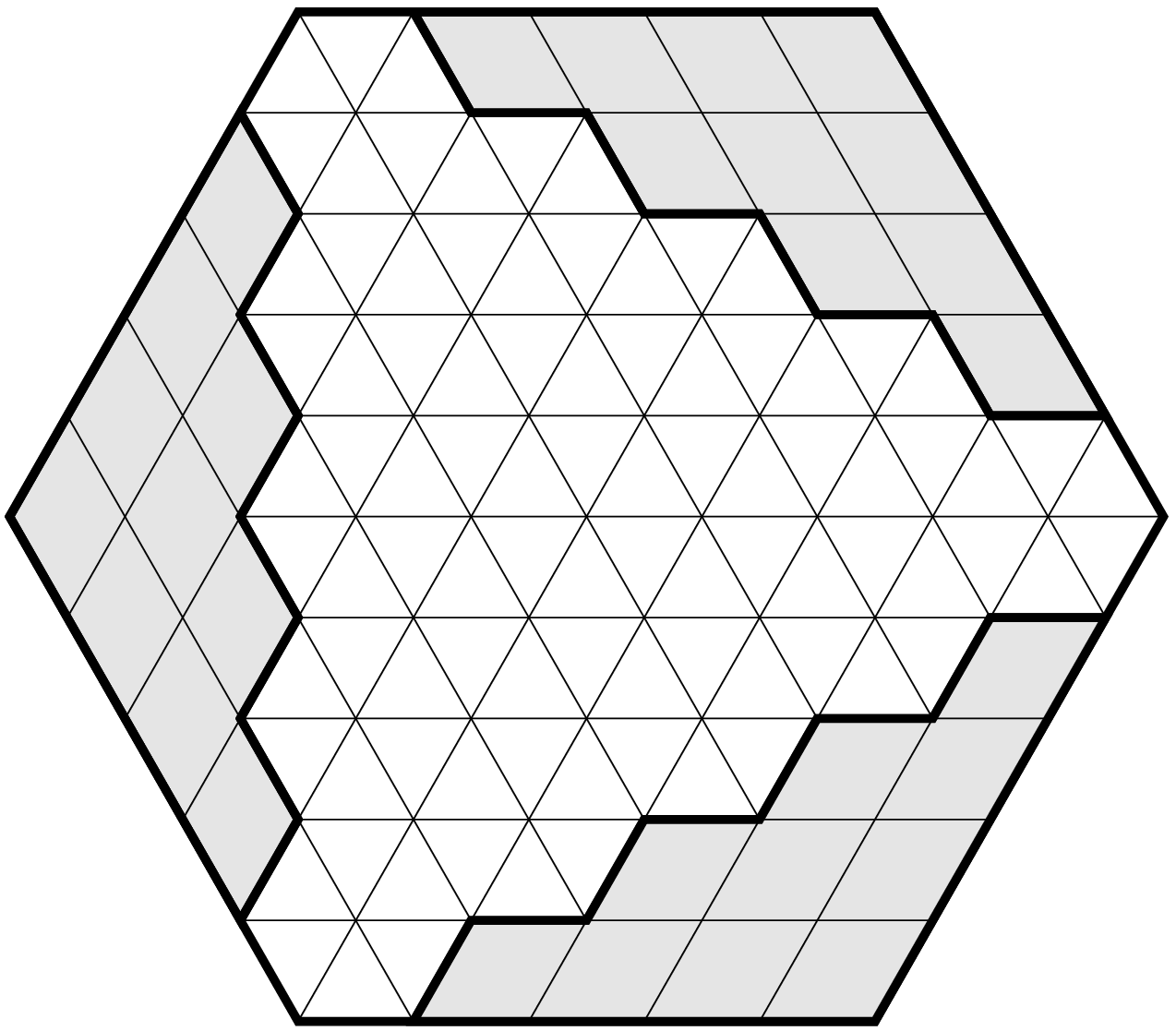}}
\twoline{\ \ \ \ \ Figure~1.4.{\rm \ \ \ \ \  }}
{\ \ \ \ \ \ \ \ \ \ \ \ \ \ \ \ Figure~1.5. {\rm $T_5$.}}
\endinsert

\proclaim{Problem 1.5} Find a formula for the number of lozenge
tilings of $T_a$%
, that explains the large amount of prime factorization of these numbers.
\endproclaim

Let us now go back to the case  
of plane partitions with the property
that their projection on two coordinate planes is contained in the
corresponding staircases. As we saw, in terms of 
tilings this
amounts to counting the tilings of a hexagon with maximal staircases removed from two
non-adjacent and non-opposite vertices. What if we remove them from adjacent vertices?
This leads to the region $H_{\text a}(a,b,c)$ (illustrated in Figure~1.6%
; evidently, the subscript stands for ``adjacent''), obtained
from $H(a,b,c)$ by removing maximal staircases from the top two corners (we assume that
these two staircases do not interfere --- otherwise the leftover region has no lozenge
tilings; non-interference amounts to $b\geq a+c-1$).

%\medskip
The $H_{\text a}(a,b,c)$'s form a family of regions that resemble the case
$b\leq a\leq c$ of the $H_{\text d}(a,b,c)$'s (compare Figures~1.6 and 1.3(c)). 
Even though they are
different, it turns out that the enumeration of tilings of both families of regions 
reduces to the evaluation of the same determinant, the one in
Theorem~1.10. We have the 
following result. 

\proclaim{Theorem 1.6} For $b\geq a+c-1$, we have
$$L(H_{\text a}(a,b,c))=P_a(b-a,c),$$
where $P_n(x,y)$ is the product on the right hand side of\/ {\rm(1.8)}.
\endproclaim

\topinsert
\centerline{\mypic{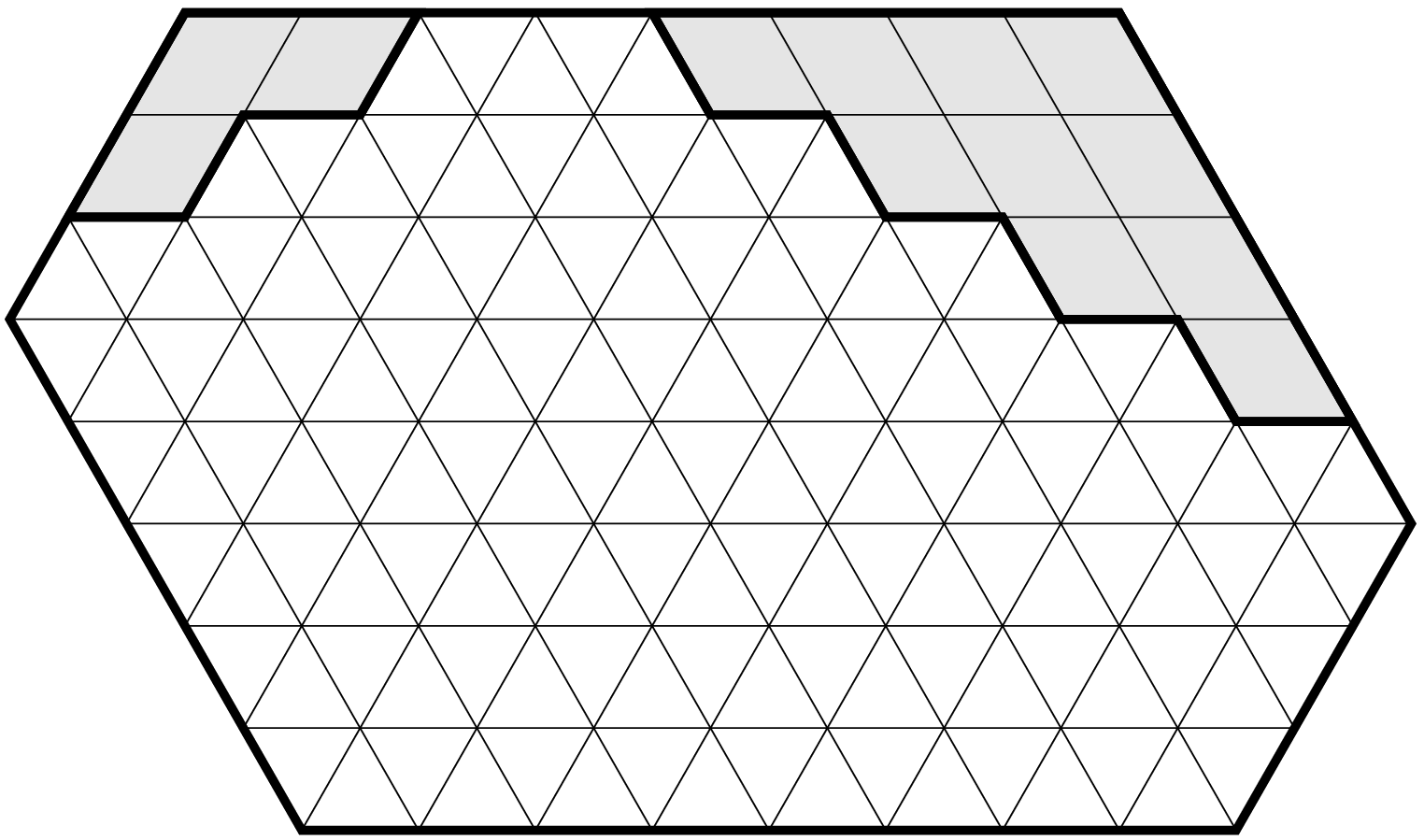}}
\centerline{Figure~1.6. {\rm $H_{\text a}(3,8,5)$.}}
\endinsert

Also here, let us consider weighted enumerators. As before, let us
denote by $L_*^*$, $L^*$, and $L_*$ the weighted tiling enumerators 
where a superscript (respectively,
subscript) indicates that the tiles along the northwestern zig-zag
line
(respectively, northeastern) are weighted by $1/2$. 
While $L^*(H_{\text a}(a,b,c))$
and $L_*(H_{\text a}(a,b,c))$ do not seem to be given by simple
product formulas, this is the case for $L_*^*(H_{\text
a}(a,b,c))$.

\proclaim{Theorem 1.7} For $b\geq a+c-1$, we have
$$\multline
L^*_*(H_{\text a}(a,b,c))=
2^{-c-a}\frac {(b+2c-a+2)_a} {(a+b-c+2)_a}\\
\times
\prod _{j=1} ^{a}\frac {(j-1)! \,(b+c-a+2j-1)! \,(b-a-c+2j+2)_j 
\,(b+2c-a+3j+2)_{a-j}} {(b+2j-1)! \,(c+a-j)!}.
\endmultline$$
\endproclaim

In view of the previously made observation that in the case that
$b<a+c-1$ the maximal staircases interfere and their removal results
in a region which is not tilable, it may seem absurd to insist on
having ``analogues" of Theorems~1.6 an 1.7 for $b<a+c-1$. But why,
instead of removing {\it maximal\/} staircases, not remove 
{\it partial\/}
staircases? To be precise, if $a+b+c\equiv 0$ mod 2, then
let us remove the partial staircase $(a-1,a-2,\dots,(a-b+c)/2)$
from the top-left vertex of the hexagon, 
and the partial staircase $(c-1,c-2,\dots,(c-b+a)/2)$ 
from the top-right vertex. 
(See Figure~1.7(a), in which the removed
staircases are indicated by the white regions. The shades should be
ignored at this point.) We obtain a region that looks like a pentagon
with an ``artificial" peak glued on top. Any lozenge tiling of
this region is uniquely determined in the rhombus that is composed
out of the triangular peak and its upside-down mirror image. 
(In Figure~1.7(a) this rhombus is shaded, and the
unique way to tile this rhombus is shown.) Hence, we may equally well
remove this rhombus. The leftover region now has the form of a
pentagon with a notch (see Figure~1.7(b); at this point the ellipses
are without relevance). Let us denote this region by $H_{\text
n}(a,b,c)$. Remarkably, extensive computer calculations suggest that
the number of lozenge tilings of $H_{\text n}(a,b,c)$ is given by a
``simple" product formula.
We state it as Conjecture~A.1 in the
Appendix. The fact that the result, 
even though given in terms of a
completely explicit product, is unusually
complex\footnote{No nontrivial simplifications seem to be possible.
We are not aware of any other ``nice" result which is similarly
involved.} may indicate that proving the conjecture may be a
formidable task.

%Christian: Relabelled figures 1.7 and 1.8, and moved the
%   new figure 1.8 to a later place.

\topinsert
\twoline{\mypic{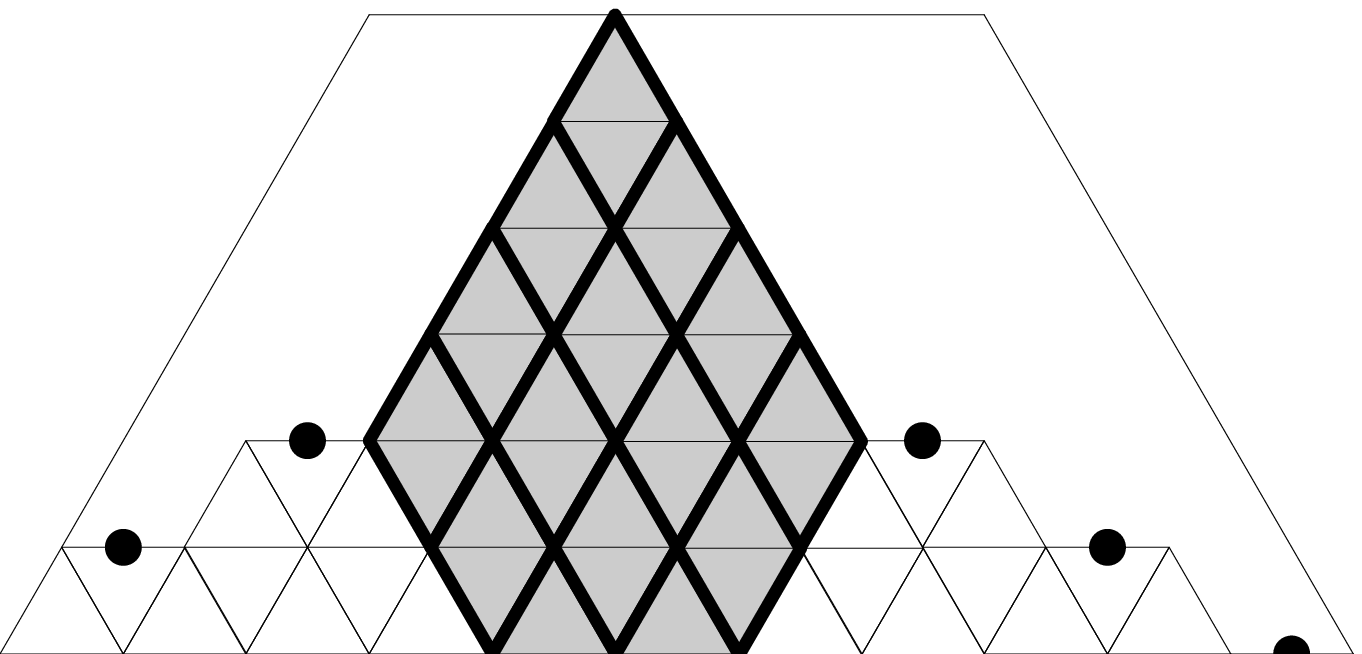}}{\mypic{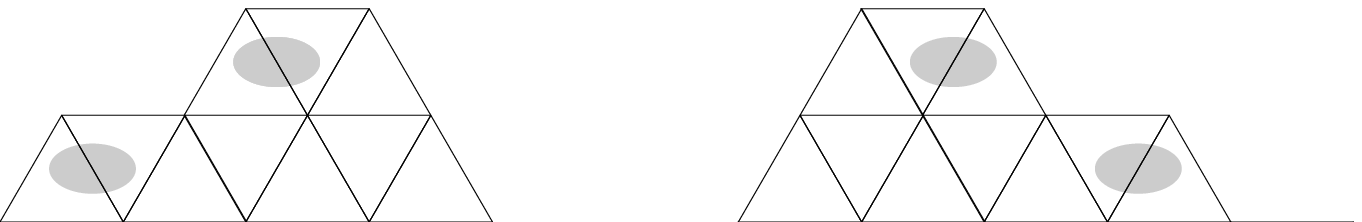}}
\twoline{Figure~1.7{\rm(a). \ \ \ \ \ \ }}{
Figure~1.7{\rm(b). $H_{\text n}(6,5,7)$}}
\vskip-5.3cm
\hbox{\hskip1.65cm
\btexdraw
  \drawdim truecm \setunitscale1.5
  \linewd.01
\htext(0.6 1){$a$}
\htext(-.75 0){$b$}
\htext(-2.1 1){$c$}
\htext(1.5 2.8){$(a+b-c)/2$}
\htext(4 2.8){$(c+b-a)/2$}
\htext(5 1){$a$}
\htext(3.6 -.1){$b$}
\htext(2 1){$c$}
\rtext td:-60 (2.1 2.6){$\left\{ \vbox{\vskip.9cm} \right. $}
\rtext td:-120 (4.7 2.5){$\left\{ \vbox{\vskip1.2cm} \right. $}
\rtext td:30 (2.2 1){$\left\{ \vbox{\vskip1.6cm} \right. $}
\rtext td:-30 (4.7 1) {$\left. \vbox{\vskip1.5cm} \right\} $}
\rtext td:90 (3.7 0.1){$\left\{ \vbox{\vskip1.2cm} \right. $}
\etexdraw
}
\vskip1cm
\endinsert

Moreover, it seems that the region $H_{\text n}(a,b,c)$ also allows a
weighted enumeration which is given by a simple product formula. Let
us, as before, denote by 
$L_*^*$, $L^*$, and $L_*$ the weighted tiling enumerators 
where a superscript (respectively,
subscript) indicates that the tiles along the northwestern zig-zag
line
(respectively, northeastern) are weighted by $1/2$. (In Figure~1.7(b)
these tiles are marked by ellipses.) While $L^*(H_{\text n}(a,b,c))$
and $L_*(H_{\text n}(a,b,c))$ do not seem to be given by simple
product formulas, this seems to be the case for $L_*^*(H_{\text
n}(a,b,c))$. We state it as Conjecture~A.2 in the Appendix. 
Again, the result is
unusually complex, which may indicate that a proof may be a
considerable undertaking.

\smallpagebreak
There is one more possibility for choosing two corners of the hexagon from which to
remove maximal staircases --- two {\it opposite} corners. Data suggests that in general
this does not lead to simple product formulas. There is one exception, when the sides
supporting the removed staircases are equal (see Figure~1.8 for an example), 
but this is a ``semi-frozen''
situation --- each tiling decomposes in tilings of parallel strips of width two (see
the proof of Proposition~1.8 in Section 2). 

\topinsert
\centerline{\mypic{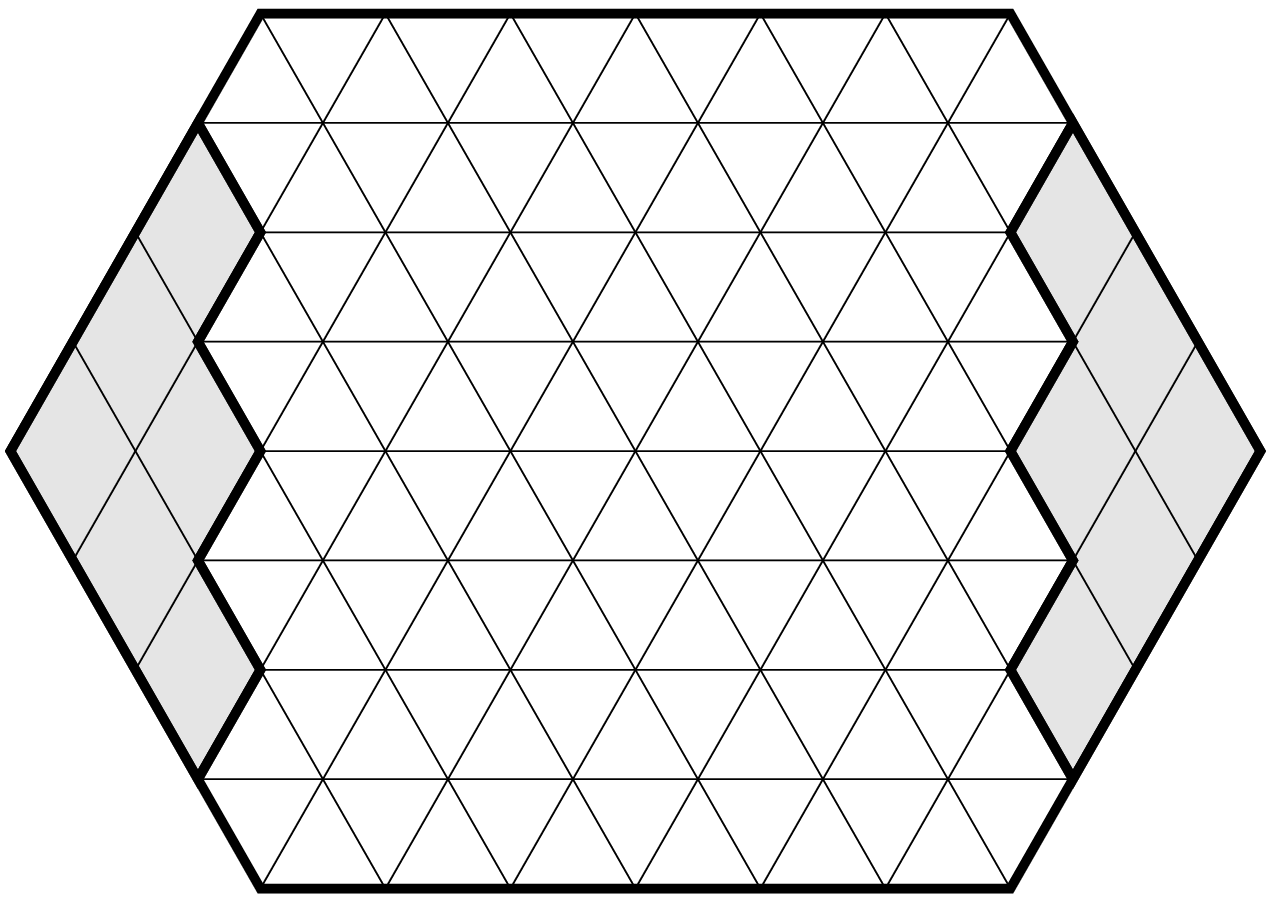}}
\centerline{Figure~1.8. {\rm $H_{\text o}(4,6,4)$.}}
\endinsert

%\medskip
Let $H_{\text o}(a,b,c)$ be the region obtained from $H(a,b,c)$ by removing maximal
staircases from the western and eastern corners. (Not unexpectedly,
here, the subscript stands for ``opposite.'')

\proclaim{Proposition 1.8}  $L(H_{\text o}(a,b,a))=(b+1)^a.$
\endproclaim

\topinsert
\twoline{\mypic{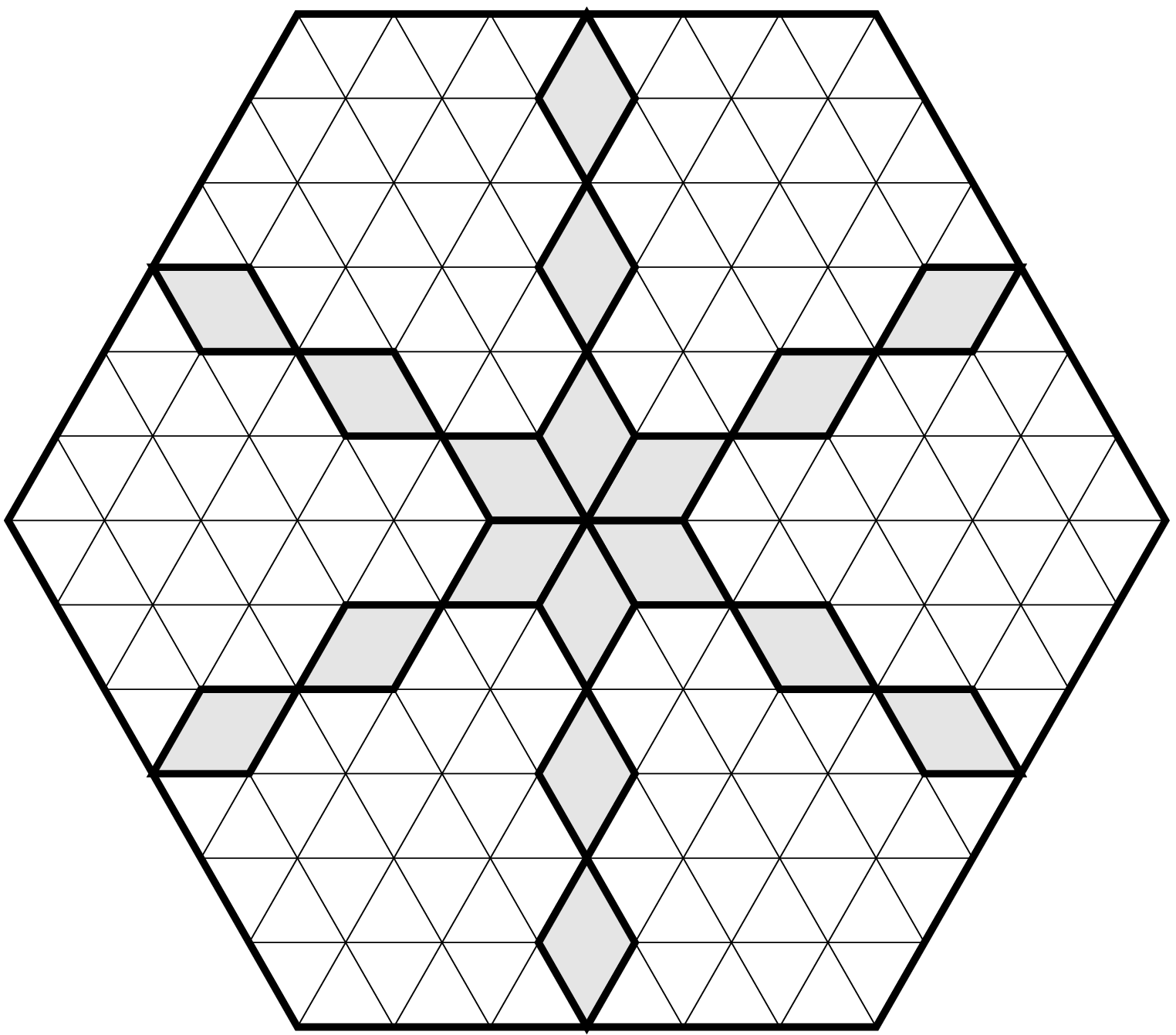}}{\mypic{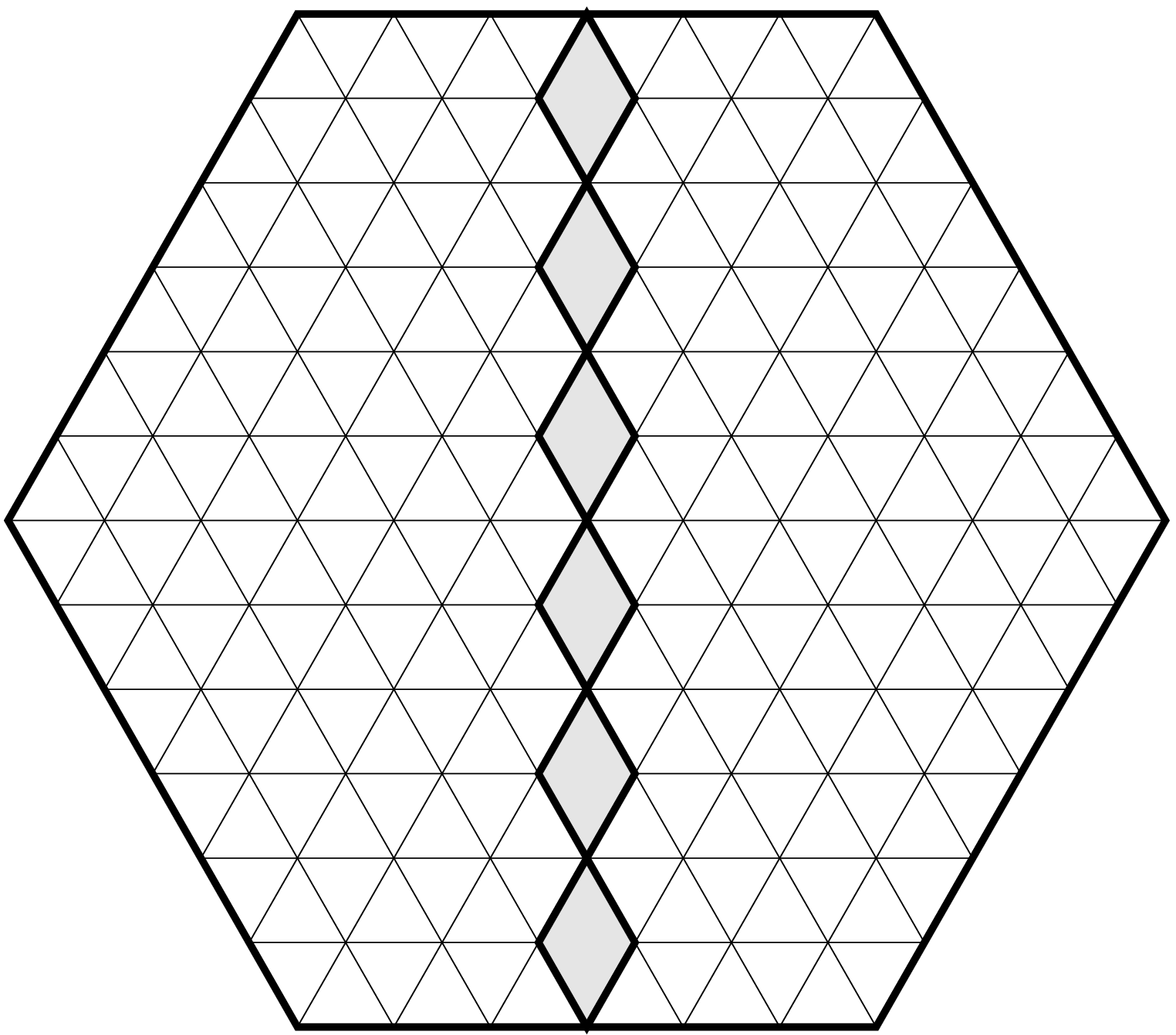}}
\twoline{Figure~1.9.{\rm \ \ \ \ \  }}{\ \ \ \ \ \ \ \ \ \
Figure~1.10.{\rm }}
\endinsert

A different viewpoint that can naturally lead one to consider the regions
introduced above (hexagons with
corners cut off) is based on symmetry classes of plane partitions. (It was in fact
this viewpoint that provided the original motivation to study these regions.) Consider
the regular hexagon $H_{2n}:=H(2n,2n,2n)$. The ten symmetry classes of plane
partitions contained in a cube of side $2n$ (see \cite{\Stapp} for their definition) are
identified with the ten symmetry classes of tilings of $H_{2n}$ (see \cite{\Ku}). 
Define a {\it ray} of tiles to be a sequence of $n$ tiles extending from the center of 
$H_{2n}$ to the nearest point of one of its edges. The six rays of tiles of $H_{2n}$ 
are shown, for $n=3$, in Figure~1.9. It is easy to see that the tilings of $H_{2n}$ 
that have CSTC 
symmetry (i.e., cyclically symmetric, transposed complementary tilings)
contain 
the tiles of all six rays, and the restriction of such a tiling to one of the six congruent regions 
left by removing the rays determines the whole tiling uniquely (see
Figure~1.9). The 
regions $H_{\text d}(a,b,c)$, $b\leq a\leq c$, form a two parameter generalization of 
this (compare Figures~1.9 and 1.3(c)).

%\medskip
Similarly, TC (i.e., transposed complementary) symmetry forces inclusion of two
opposite rays in the tiling, and reduces to enumerating tilings of one of the two
pieces left over after removing the two opposite rays of tiles (see
Figure~1.10). 
The regions
corresponding to the plane partitions considered by Proctor (hexagons with one corner
cut off) form a one parameter generalization of this (compare
Figures~1.10 and 1.1).

\topinsert
\twoline{\mypic{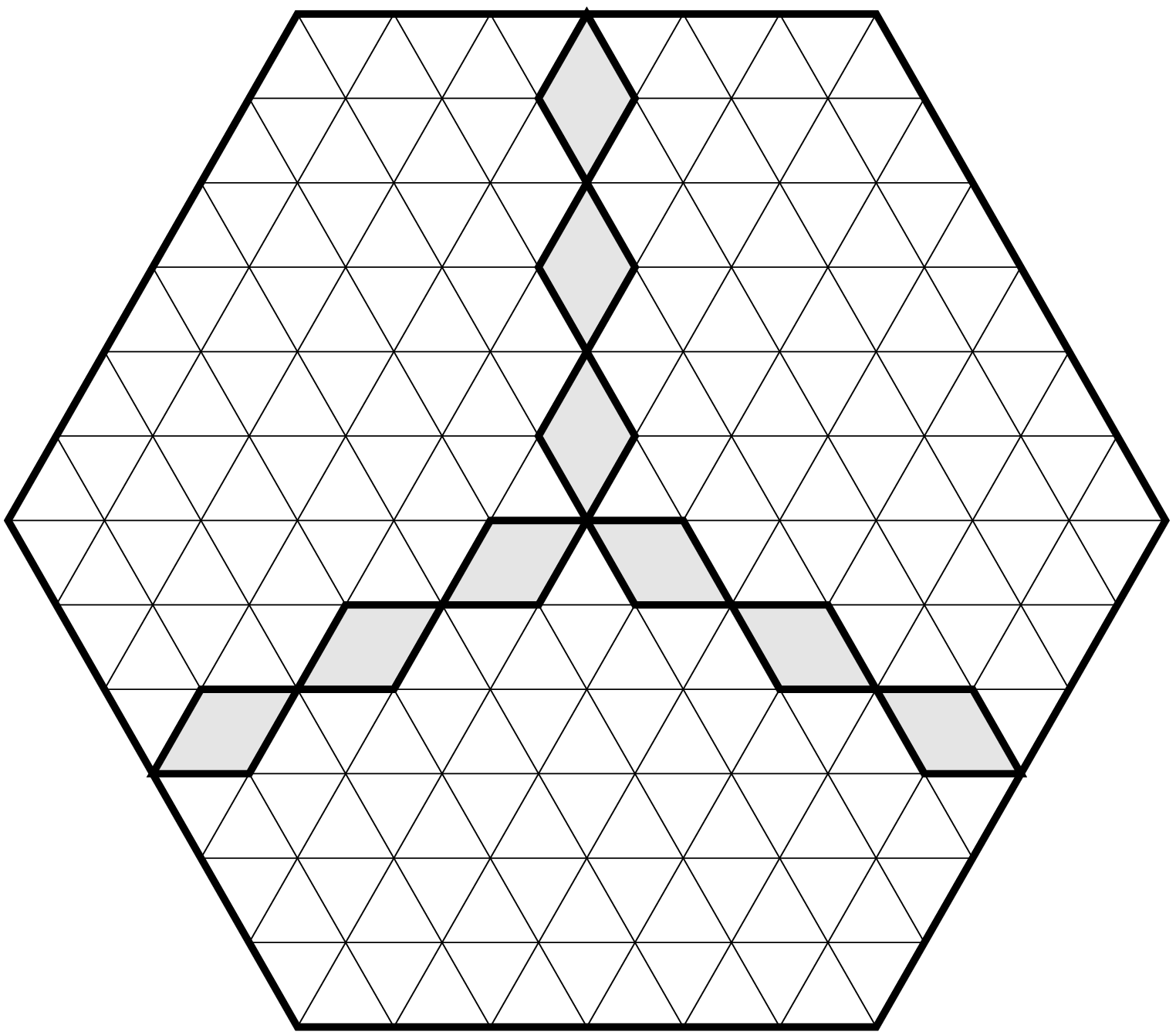}}{\mypic{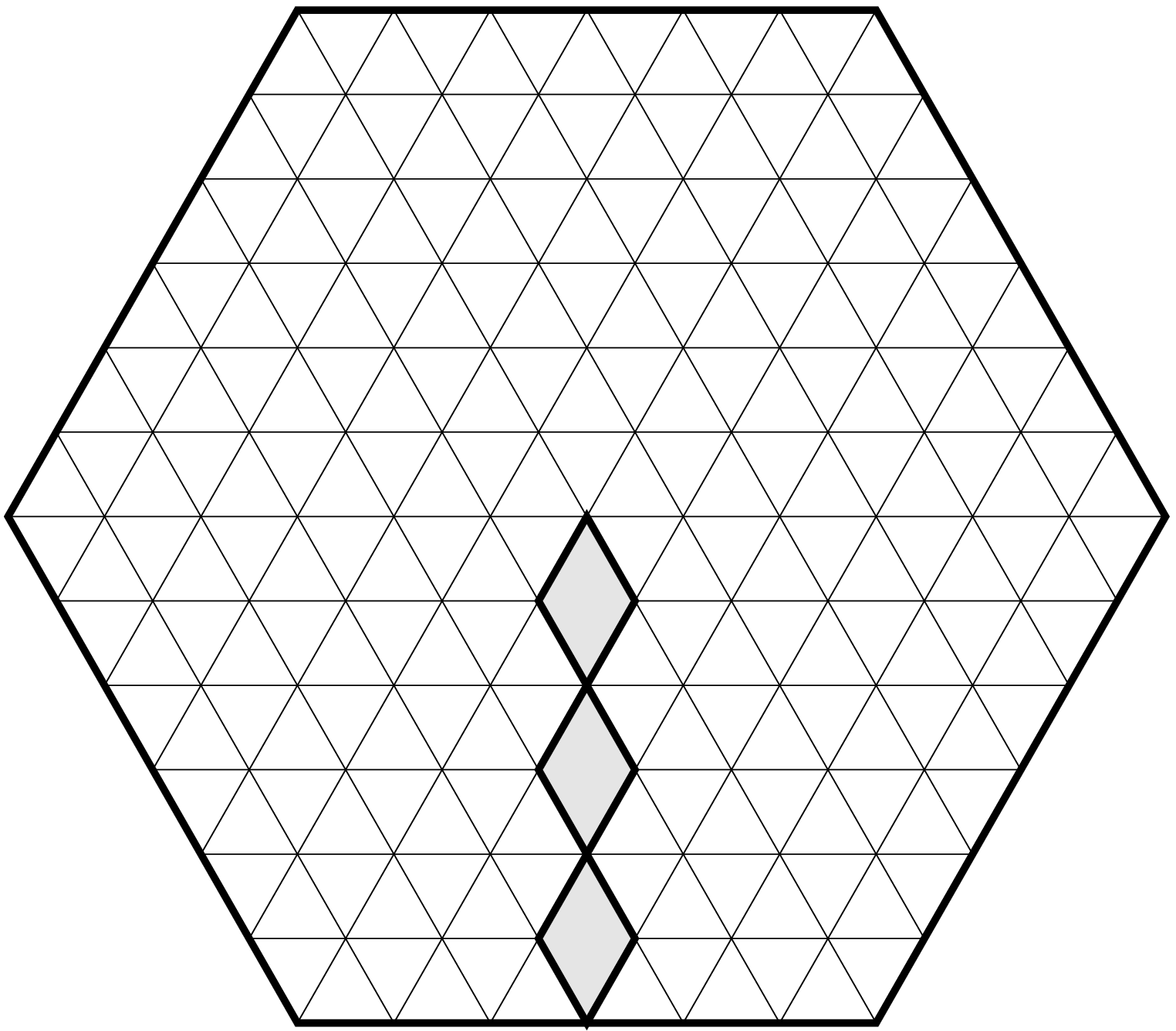}}
\twoline{Figure~1.11.{\rm \ \ \ \ \  }}{\ \ \ \ \ \ \ \ \ \
Figure~1.12.{\rm }}
\endinsert

We are thus naturally led
to consider the regions generated by removing {\it three}
alternating rays, as shown in Figure~1.11. This region does not correspond to a symmetry
class of plane partitions, but it is nevertheless quite 
compelling
to consider in this 
context. The regions $H_{\text a}(a,b,c)$ form a two parameter generalization of it
(compare Figures~1.11 and 1.6).
The two parameters were essential in conjecturing a formula for the number of tilings
of these regions, based on data: polynomials fully factored into linear factors
contain much more information than just integers factored into small primes. 

%\medskip
To finish this analysis, we mention that removing {\it one} ray (see
Figure~1.12) 
leads to a region
whose number of tilings has a simple product expression, as it easily follows from the
Factorization Theorem of \cite{\Ci1} and the results of \cite{\Cipp1}. 

%\medskip
To rephrase the above statements, if one regards the
regions formed by removing rays as being built up of $60^\circ$ sectors, the
tilings of the regions consisting of 1, 2, 3 or 6 sectors are enumerated by simple
product formulas. Data suggests that this is not the case for 4 or 5
sectors.

We prove our results in Section~2 
by employing a standard bijection between lozenge
tilings and 
non-intersecting lattice paths, thus, due to the
Lindstr\"om--Gessel--Viennot theorem \cite{\LindAA}\cite{\GV},
obtaining a determinant for the number (respectively weighted count) 
of tilings that we are
interested in, and, finally, by evaluating the resulting determinants.
In most of the cases, we obtain special cases of the two 
determinant evaluations that we state in Theorems~1.9 and 1.10 below.
Theorem~1.9 is due to the second author \cite{\KratBD, (5.3)} (see 
\cite{\Amdeb} for a simple proof, which is reproduced in \cite{\CKpp2}). 
The determinant evaluation in
Theorem~1.10 
does not seem to have appeared previously in the literature.
The paper
\cite{\CiKrAD} contains our original proof, which is rather involved,
but has its own appeal as it contains a non-automatic (!) application
of Gosper's algorithm \cite{\GospAB} (see also 
\cite{\GrKPAA, Sec.~5.7}\cite{\PeWZAA, Sec.~II.5}). Later we
discovered that, in fact, there is a combinatorial argument which
transforms the determinant in Theorem~1.10 into an instance of the
determinant in Theorem~1.9, so that these two determinant evaluations
are actually equivalent. It is this argument that we give in
Section~2.

\proclaim{Theorem 1.9}Let $x$, $y$ and $n$ be nonnegative integers with
$x+y>0$. Let $K_n(x,y)$ be the matrix 
$$\align
K_n(x,y)&:=
\left(\frac{(x+y+i+j-1)!}{(x+2i-j)!\,(y+2j-i)!}\right)_{1\leq i,j\leq
n}\\
&\hphantom{:}=
\left(\frac{1}{x+2i-j}\binom {x+y+i+j-1}{y+2j-i}\right)_{1\leq i,j\leq
n}.
\tag1.5\endalign$$
Then 
$$\det(K_n(x,y))=\prod_{j=1}^{n}\frac{(j-1)!\,(x+y+j)!\,
(2x+y+2j+1)_{j-1}\,(x+2y+2j+1)_{j-1}}
{(x+2j-1)!\,(y+2j-1)!}.\tag1.6$$
where, as before, the shifted factorial
%Christian: To be consistent with the other places, I changed
%  the a's here into \alpha's.
$(\alpha)_k$ is defined by $(\alpha)_k:=\alpha(\alpha+1)\cdots
(\alpha+k-1)$, $k\ge1$, and $(\alpha)_0:=1$.
\endproclaim

\proclaim{Theorem 1.10} Let $n$ be a positive integer, and let $x$ and
$y$ be nonnegative integers. Let $A_n(x,y)$ be the matrix 
$$A_n(x,y):=\left(\binom{x+y+j}{x-i+2j}-\binom{x+y+j}{x+i+2j}\right)_{1\leq i,j\leq n}.
\tag1.7$$
Then 
$$\det A_n(x,y)=\prod _{j=1} ^{n}\frac {(j-1)!\,(x+y+2j)!\,(x-y+2j+1)_j\,
(x+2y+3j+1)_{n-j}} {(x+n+2j)!\,(y+n-j)!}.\tag1.8$$
\endproclaim

\flushpar
{\bf Note.} 
Theorems~1.9 and 1.10 are only formulated for
nonnegative integral $x$ and
$y$. But in fact, with a
generalized definition of factorials and binomials (cf\.
\cite{\GrKPAA, \S5.5, (5.96), (5.100)},
both theorems would also make sense and be true for complex $x$ and $y$.

\mysec{2. Proofs of the results}

\medskip

As already mentioned at the end of the Introduction,
we employ in our proofs a standard bijection that maps each lozenge tiling $T$ of
a region $R$ on the triangular lattice to a family of non-intersecting lattice paths 
taking steps east or north on the grid lattice $\Z^2$. This bijection works as
follows. Choose a lattice line direction $d$ --- without loss of generality, the
horizontal direction --- , and mark with a dot the unit segments parallel to $d$ on the 
boundary of $R$ (see Figure~2.1). 
Each marked segment is either on top or at the bottom of the region $R$.
Label from right to left the ones at the bottom of $R$ by
$u_1,u_2,\dotsc,u_m$ and the
ones on top of $R$ by $v_1,v_2,\dotsc,v_n$.

\topinsert
\vbox{
\twoline{\mypic{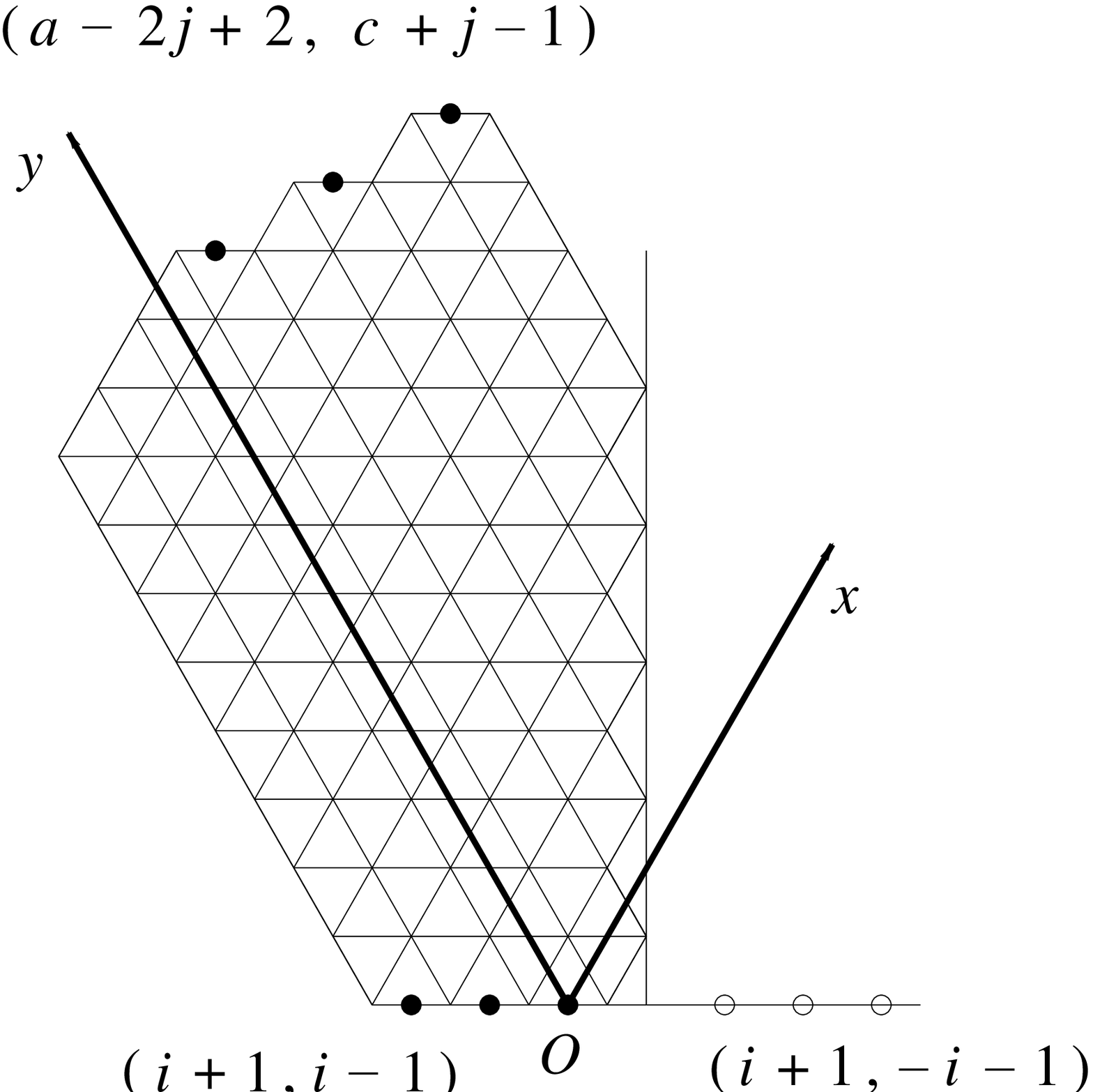}}{\mypic{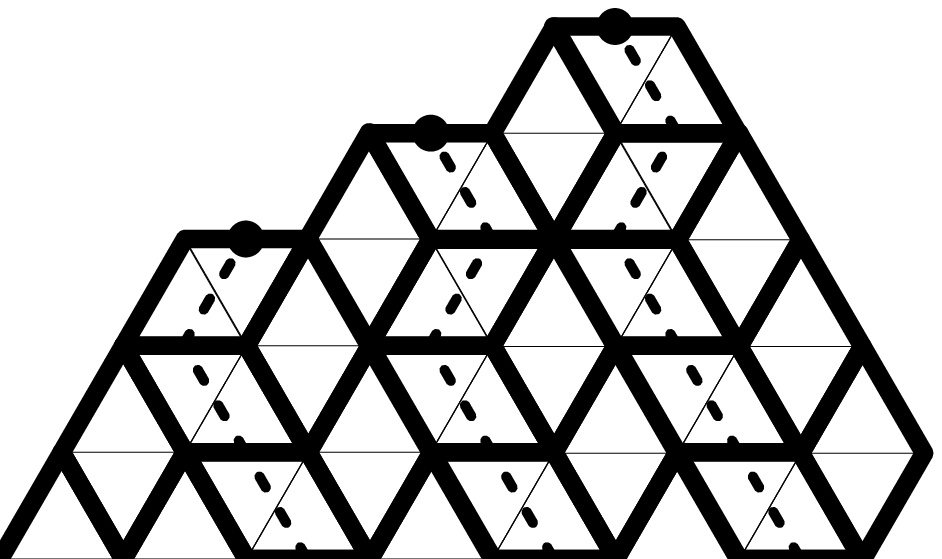}}
\twoline{Figure~2.1.{\rm \ \ \ \ \  }}{\ \ \ \ \ \ \ \ \ \ Figure~2.2.{\rm }}
\vskip-1.2cm
\centerline{\hskip7.2cm$u_3\ u_2\ u_1$}
\vskip-5.9cm
\centerline{\hskip6.7cm$v_1$}
\vskip-0.05cm
\centerline{\hskip5.4cm$v_2$}
\vskip-0.05cm
\centerline{\hskip4.1cm$v_3$}
\vskip5cm
}
\endinsert

\topinsert
\vskip10pt
$$
\Einheit=.4cm
\Gitter(7,11)(-3,-2)
\Koordinatenachsen(7,11)(-3,-2)
\Pfad(0,0),1212211222212\endPfad
\Pfad(-1,1),212212212212\endPfad
\Pfad(-2,2),22221221221\endPfad
\Diagonale(-1,-3){7}
\DickPunkt(0,0)
\DickPunkt(-1,1)
\DickPunkt(-2,2)
\DickPunkt(5,8)
\DickPunkt(3,9)
\DickPunkt(1,10)
\Label\ru{y=x-2}(3,0)
\Label\lu{u_1}(0,0)
\Label\lu{u_2}(-1,1)
\Label\lu{u_3}(-2,2)
\Label\ro{v_1}(5,8)
\Label\ro{v_2}(3,9)
\Label\ro{v_3}(1,10)
\hbox{\hskip5.5cm}
\Gitter(7,11)(-3,-3)
\Koordinatenachsen(7,11)(-3,-3)
\Pfad(0,-1),21212211222212\endPfad
\Pfad(-1,-2),222212212212212\endPfad
\Pfad(-2,-3),2222222221221221\endPfad
\Diagonale(-2,-4){8}
\DickPunkt(0,0)
\DickPunkt(-1,1)
\DickPunkt(-2,2)
\DickPunkt(0,-1)
\DickPunkt(-1,-2)
\DickPunkt(-2,-3)
\DickPunkt(5,8)
\DickPunkt(3,9)
\DickPunkt(1,10)
\Label\ru{y=x-2}(3,0)
\hskip1.8cm
$$
\vskip7pt
\twoline{Figure~2.3{\rm(a). \ \ \ \ \  }}{\ \ \ \ \ \ \ \ \ \
Figure~2.3{\rm(b).  }}
\vskip10pt
\endinsert

Consider now the tile $t_1$ of $T$ resting on $u_i$, for some fixed $1\leq i\leq m$
(see Figure~2.2 for an example of a tiling $T$ of the region from Figure 2.1).
Let $t_2$ be the other tile of $T$ containing the side of $t_1$ opposite $u_i$.
Continue the sequence of selected tiles by choosing $t_3$ to be the tile of $T$ sharing
with $t_2$ the side of $t_2$ opposite the one common to $t_2$ and $t_1$. This
procedure leads to a path of rhombic tiles growing upward, which clearly must end on
one of the $v_j$'s. 
This path of rhombi can be identified with a (linear) path that
starts at the midpoint of $u_i$ and ends at the midpoint of $v_j$,
see Figure~2.2. (There, the resulting linear paths are indicated by
dotted segments.) After normalizing the oblique coordinate system and
rotating it in standard position, we obtain a lattice path on
$\Z^2$ that starts at the midpoint of $u_i$ (actually, its image
after these normalizations), ends at the midpoint of $v_j$ (again, actually 
its image after these normalizations) and takes 
unit steps east or north. (See Figure~2.3(a) for the resulting
lattice paths in our example.)
One obtains this way a family $\Cal P$ of $m$ 
lattice paths, one for each $1\leq
i\leq m$, and they cannot touch each other
since the corresponding paths of tiles are disjoint.
We obtain in particular that $m\leq n$, and hence by symmetry $m=n$. It is easy to
see that the correspondence $T\mapsto \Cal P$ is a bijection between the set of
tilings $T$ of $R$ and the families $\Cal P$ 
of non-intersecting lattice paths starting at the
midpoints of $u_1,u_2,\dotsc,u_n$, ending at the midpoints of
$v_1,v_2,\dotsc,v_n$ and
contained within $R$.

\proclaim{Lemma 2.1} For $b\leq a\leq c$, we have 
$L(H_{\text d}(a,b,c))=(-1)^{b(b+1)/2} \det A_b(a-2b-1,b+c+1)$.
\endproclaim

\pf Label the horizontal unit segments on the boundary of $H_{\text d}(a,b,c)$ as
described above. 
%This yields segments $u_1,u_2,\dotsc,u_b$ below our region and segments
%$v_1,v_2,\dotsc,v_b$ above it.
Choose an oblique coordinate system with the origin at the midpoint of $u_1$ and
coordinate axes parallel to the non-horizontal lattice lines of the triangular
lattice (see Figure~2.1). 
Applying the procedure described above to $H_{\text d}(a,b,c)$ one obtains that 
$L(H_{\text d}(a,b,c))$ is equal to the number of families of non-intersecting
lattice paths 
with starting points $u_i=(-i+1,i-1)$, $i=1,2,\dotsc,b$, ending points
$v_j=(a-2j+2,c+j-1)$, $j=1,2,\dotsc,b$, and with the additional requirement that they 
do not touch the
line 
$y=x-2$.
(This requirement ensures that the corresponding paths of rhombi stay
within our region $H_{\text d}(a,b,c)$; see Figure~2.3(a). 
Note that by abuse of notation we denote the
midpoints of the $u_i$'s and $v_j$'s by the same symbols we use for the
segments.)

By the Lindstr\"om-Gessel-Viennot theorem 
\cite{\LindAA, Lemma~1}\cite{\GV}\cite{\Stenlp, Theorem~1.2}, the number of 
such families
of non-intersecting lattice paths is given by the determinant of the $b\times b$ matrix
whose $(i,j)$-entry is the number of lattice paths from $u_i$ to $v_j$ that
are strictly above the line $y=x-2$. By Andr\'e's Reflection Principle \cite{\Feller}, 
this number is
$$\binom{a+c-j+1}{a-2j+i+1}-\binom{a+c-j+1}{a-2j-i+1}.$$
Therefore, the Lindstr\"om-Gessel-Viennot theorem implies that
$$L(H_{\text d}(a,b,c))=
\det \left(\binom{a+c-j+1}{a-2j+i+1}-\binom{a+c-j+1}{a-2j-i+1}
\right)_{1\leq i,j\leq b}.\tag2.1$$
Reversing columns (i.e., replacing $j$ by $b+1-j$) in (2.1) the right
hand side becomes
$$(-1)^{b(b-1)/2}\det \left(\binom{a+c+j-b}{a+2j+i-2b-1}-\binom{a+c+j-b}{a+2j-i-2b-1}
\right)_{1\leq i,j\leq b}.\tag2.2$$
The entries of the matrix in (2.2) are readily seen to be precisely the negatives of
the entries of $A_b(a-2b-1,b+c+1)$. This implies the statement of the Lemma.  $\square$

%$$\align
%(-1)^{b(b-1)/2}\det \left(\binom{a+c+j-b}{a+2j+i-2b-1}-\binom{a+c+j-b}{a+2j-i-2b-1}
%\right)_{1\leq i,j\leq b}&\\
%(-1)^{b(b-1)/2}\det
%\left(\binom{(a-2b-1)+(b+c+1)+j}{(a-2b-1)+2j+i}-\binom{(a-2b-1)+(b+c+1)+j}
%{(a-2b-1)+2j-i}
%\right)_{1\leq i,j\leq b}&\\
%(-1)^{b+b(b-1)/2}\det A_b(a-2b-1,b+c+1)&,
%\endalign$$

\smallpagebreak
{\it Proof of Theorem 1.10.} The preceding proof of Lemma~2.1 shows
(by renumbering the lattice paths from left to right)
that $(-1)^{n(n+1)/2}A_n(x,y)$ counts the number of families $\Cal P$
of $n$ non-intersecting lattice paths, with starting points
$(-n+i,n-i)$, $i=1,2,\dots,n$, and end points $(x+2j+1,y-j-1)$,
$j=1,2,\dots,n$,
with the additional requirement that they do not touch the line
$y=x-2$ (see Figure~2.3(a)
for an example). Now, for each such family, 
we prepend $(2n-2i+1)$ vertical steps to the
$i$-th path. Thus we obtain families $\Cal P'$
of $n$ non-intersecting lattice paths, with starting points
$(-n+i,-n+i-1)$, $i=1,2,\dots,n$, and end points 
$(x+2j+1,y-j-1)$, $j=1,2,\dots,n$,
with the additional requirement that they do not touch the line
$y=x-2$. (See Figure~2.3(b)
for the resulting path family in our example.) In fact, this is a
bijection between the former and latter path families, because the
prepended path portions are in fact forced by the boundary $y=x-2$. 
Therefore $(-1)^{n(n+1)/2}A_n(x,y)$ 
is also equal to the number of the latter path
families. 

Again we may apply the Lindstr\"om-Gessel-Viennot theorem. Since the
Reflection Principle yields that the number of paths from
$(-n+i,-n+i-1)$ to $(x+2j+1,y-j-1)$ which do not touch the line
$y=x-2$ is given by
$$\multline
\binom {x+y+j+2n-2i+1}{x+2j+n-i+1}-\binom {x+y+j+2n-2i+1}{x+2j+n-i}\\
=
\frac {(y-x-3j)\,(x+y+j+2n-2i+1)!} {(x+2j+n-i+1)!\,(y+n-j-i+1)!},
\endmultline$$
we infer that $(-1)^{n(n+1)/2}A_n(x,y)$ is equal to
$$\align
\det\bigg(&\frac {(y-x-3j)\,(x+y+j+2n-2i+1)!}
{(x+2j+n-i+1)!\,(y+n-j-i+1)!}\bigg)_{1\le i,j\le n}\\
&=\prod _{j=1} ^{n}(y-x-3j)\cdot
\det\left(\frac {1} {y+n-j-i+1}\binom {x+y+j+2n-2i+1}
{x+2j+n-i+1}\right)_{1\le i,j\le n}\\
&=\prod _{j=1} ^{n}(y-x-3j)\cdot
\det\left(\frac {(-1)^{x+2j+n-i+1}} {y+n-j-i+1}\binom {-y-n+j+i-1}
{x+2j+n-i+1}\right)_{1\le i,j\le n}\\
&=\prod _{j=1} ^{n}\big((-1)^{x+n-j}(y-x-3j)\big)\\
&\hskip2cm\times
\det\left(\frac {1}
{-y-x-2n+2i-j-2}\binom {-y-n+j+i-2}{x+2j+n-i+1}\right)_{1\le i,j\le n}.
\endalign$$
(At the second to last equality we used that
$\binom{n}{k}=(-1)^k\binom{-n+k-1}{k}$.)
This latter determinant is the determinant $\det K_n(-y-x-2n-2,x+n+1)$, 
with $K_n(x,y)$ defined in (1.5). 
In view of Theorem~1.9, this proves (1.8), after
some manipulation. $\square$

\smallpagebreak
{\it Proof of Theorem 1.1.} This follows directly from Lemma~2.1 and Theorem~1.10. 
$\square$

\smallpagebreak
{\it Proof of Theorem 1.2.} As in the proof of Lemma~2.1,
we map the tilings to families $\Cal P$
of non-intersecting lattice paths, with starting points 
$u_i=(-i+1,i-1)$, $i=1,2,\dotsc,b$, ending points
$v_j=(a-2j+2,c+j-1)$, $j=1,2,\dotsc,b$, and with the additional requirement 
that they do not touch the line $y=x-2$ (see again Figure~2.3(a)). 
Following the proof of
Theorem~1.10, we prepend $2i-1$ vertical steps to the $i$-th path
(because here we kept the numbering of the paths from right to left), so
that we obtain families $\Cal P'$ 
of non-intersecting lattice paths, with starting points 
$u_i=(-i+1,-i)$, $i=1,2,\dotsc,b$, ending points
$v_j=(a-2j+2,c+j-1)$, $j=1,2,\dotsc,b$, and with the additional requirement 
that they do not touch the line $y=x-2$ (see again Figure~2.3(b)). 
However, unlike in the
proofs of Lemma~2.1 and Theorem~1.10, here each family has a certain
weight, given by the $\ell$-th power of $1/2$, where $\ell$ is the number
of lozenges that are weighted by $1/2$ in the corresponding tiling
(compare Figures~1.3(c), 2.1 and 2.3(b)).

To realize this weight, we give each horizontal step from
$(a-2k+1,c+k-1)$ to $(a-2k+2,c+k-1)$, for some $k$, 
the weight $1/2$. This takes
care of the fact that the lozenges along the northwestern zig-zag
line are weighted by $1/2$. To also take into account that the
lozenges along the eastern zig-zag line are weighted by $1/2$, we
assign a weight of 2
%Christian:
(sic!) to each horizontal step from
%Christian: Changed -k into k.
$(k,k)$ to $(k+1,k)$,
for some $k$
(i.e., to each horizontal step which terminates directly at the line
$y=x-1$
the paths are not allowed to cross). This generates our weight, up
to a multiplicative constant of $2^a$. 
Indeed, for each marked lozenge position along the eastern boundary, 
a tiling $T$ has
a lozenge in that position if and only if the corresponding unit 
segment weighted by 2
%Christian: modified that place.
is {\it not\/} a step on some lattice path 
of the family corresponding to $T$, thereby giving rise to a missing
weight of $1/2$ in comparison to path families where some path does
contain that step. To give an explicit example,
the tiling in Figure~2.2 contains two tiles weighted by $1/2$ 
along the eastern boundary, the third and the fifth from the
bottom. These are the ones which are completely white, as there
is no lattice path running through them. Hence, in the corresponding
path family shown in Figure~2.3(b), the step from $(2,2)$ to $(3,2)$
and the step from $(4,4)$ to $(5,4)$ (both weighted by 2) 
are not contained by 
any of the paths.

Now we want to write down the Lindstr\"om--Gessel--Viennot determinant
for our weighted count (as defined just above). 
In order to do so, we need the weighted count
of paths from $(-i+1,-i)$ to $(a-2j+2,c+j-1)$ which do not touch
$y=x-2$.
We claim that the weighted count of these latter paths is the same as the
weighted count of {\it all\/} paths from $(-i+1,-i)$ to 
$(a-2j+2,c+j-1)$, in which the last step of the
path has weight $1/2$ if it is a horizontal step. (It should be noted
that in this weighted count there are no steps with weight 2 anymore.) 
This is seen as follows: suppose we are
given a path from $(-i+1,-i)$ to $(a-2j+2,c+j-1)$ 
which does not touch the line $y=x-2$, and which has exactly $\ell$
touching points on the line $y=x-1$. These $\ell$ touching points on
$y=x-1$ must be reached by horizontal steps from
%Christian: Changed -k into k.
$(k,k)$ to $(k+1,k)$, each of which contribute a weight of $2$.
Thus, in total, this gives contribution of $2^\ell$ to the weight.
Now we map such a path to $2^\ell$ paths from
$(-i+1,-i)$ to $(a-2j+2,c+j-1)$ (without restriction), 
by focussing on the path
portions between two consecutive touching points, including the
portion between $(-i+1,-i)$ and the first touching point, and keeping
any of them either fixed or reflecting it in the line 
$y=x-1$. This proves the assertion. 

%Christian: More text.
By distinguishing between the cases where the last step of a path is
vertical respectively horizontal, 
the new weighted count of the paths from $(-i+1,-i)$ to $(a-2j+2,c+j-1)$ 
then is seen to be
$$\binom {a+c+2i-j-1}{a-2j+i+1}+\frac {1} {2}\binom
{a+c+2i-j-1}{a-2j+i}=
\frac {1} {2}\frac {(a+2c+3i-1)\,(a+c+2i-j-1)!}
{(a-2j+i+1)!\,(c+i+j-1)!}.$$

Therefore the Lindstr\"om--Gessel--Viennot determinant 
is seen, by manipulations similar to those in the proof of Theorem 1.10, to be
$$\multline
\det\bigg(\frac {1} {2}\frac {(a+2c+3i-1)\,(a+c+2i-j-1)!}
{(a-2j+i+1)!\,(c+i+j-1)!}\bigg)_{1\le i,j\le b}\\
=2^{-b}\prod _{i=1} ^{b}\big((-1)^{a+i}(a+2c+3i-1)\big)
\det\left(\frac {1} {-a-c+j-2i}\binom {-c-i-j}{a-2j+i+1}\right)_{1\le i,j\le b}.
\endmultline$$
After reversing the order of rows and columns (i.e., after replacing $i$
by $b+1-i$ and $j$ by $b+1-j$), it is seen that this
determinant is the determinant $\det K_b(-a-c-b-1,a-b)$. Application
of Theorem~1.9, division of the resulting expression by $2^a$, 
and some rearrangement finish the proof of the
theorem. $\square$

\smallpagebreak
{\it Proof of Proposition 1.3.} Consider the region $H_{\text d}(a,b,c)$ and choose the
direction $d$ in the bijection between tilings and lattice paths to be the direction
of the lattice lines going from southwest to northeast
(see Figure~2.4). Consider the unit segments 
parallel to $d$ on the boundary, and label the midpoints of those on the eastern 
boundary, from top to bottom, by $u_1,u_2,\dotsc,u_a$, and the midpoints of those on the
northwestern boundary by $v_1,v_2,\dotsc,v_a$. Choose an oblique coordinate system
centered $\sqrt{3}$ units above $u_1$ (see Figure~2.4) and with axes along the 
northwestern and western lattice line
directions.

\topinsert
\centerline{\mypic{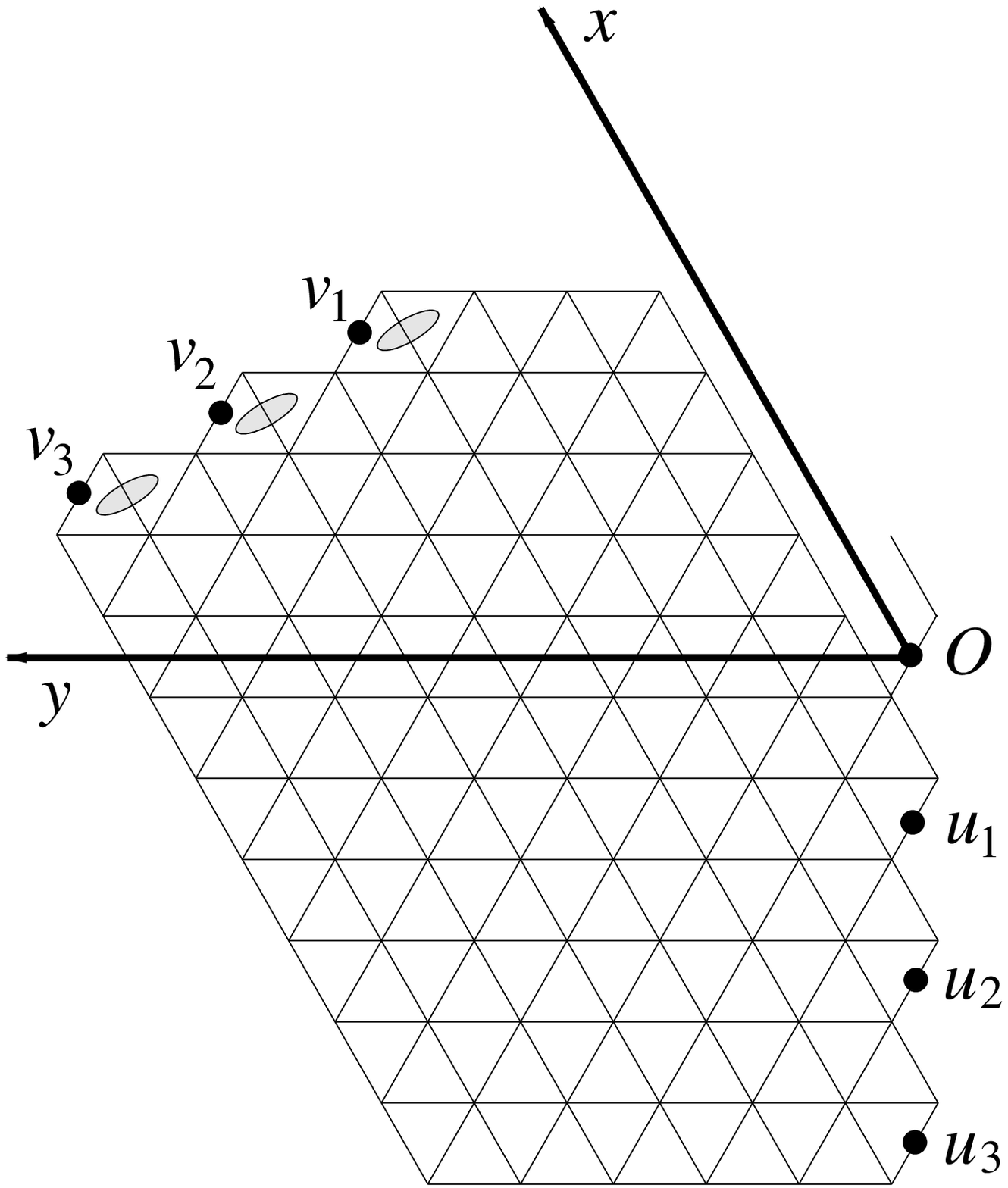}}
\centerline{Figure~2.4. {\rm }}
\endinsert

By the bijection between tilings and lattice paths, each tiling $T$ of $H_{\text
d}(a,b,c)$ is identified with a family $\Cal P$ of 
non-intersecting lattice paths with 
starting points 
%Christian: Unfortunately, the given coordinatization does not match
%  with Figure 2.4. Either you have to shift the cordinate axes
%  by two units up, or you have to add (2,-1) to the starting
%  and end points. (I consider the latter solution less attractive.) 
$(-2i,i)$, $i=1,2,\dotsc,a$, 
and ending points $(c-a-j,b-a+2j)$, 
$j=1,2,\dotsc,a$.
(Unlike in the proof of Lemma~2.1, no additional requirements on the paths are
needed). 

The weighted enumerator $L^*$ assumes that the northwesternmost $a$ tile positions
are weighted by $1/2$ (these positions are marked in Figure~2.4). 
Correspondingly, the paths of $\Cal P$ whose last step is a
%Christian: just to be consistent with the rest
%$y$-step 
vertical step have weight $1/2$ 
%Christian: minor addition
(as before, the weight of a lattice path is the product of the weights
of its steps). The weight of the family $\Cal P$ is the product of the weights of its
members, and by construction it matches the weight of the tiling $T$. 
By the Lindstr\"om-Gessel-Viennot theorem,
the total weight of those families $\Cal P$ that are
non-intersecting --- and hence $L^*(H_{\text d}(a,b,c))$ --- is given by the determinant of
the $a\times a$ matrix whose $(i,j)$-entry equals the total weight of the paths from
$u_i$ to $v_j$. Splitting the latter family in two according to the type of the last
step, one obtains that its total weight is
$$\align
&\binom{b+c-2a+i+j-1}{c-a+2i-j-1}+\frac{1}{2}\binom{b+c-2a+i+j-1}{c-a+2i-j}\\
&\ \ \ \ =\frac{b+2c-3a+3i}{2}\frac{(b+c-2a+i+j-1)\,!}
{(c-a+2i-j)\,!\,(b-a+2j-i)\,!}.
\endalign
$$
When computing the determinant of the above matrix, the $j$-free factors can be
factored out along rows, and the leftover matrix is precisely the one in
(1.5) with
$n=a$, $x=c-a$, and $y=b-a$. 
Using (1.6) and substituting the values of $x$ and $y$ one obtains the formula 
in the statement of Proposition~1.3. $\square$

\smallpagebreak

%Christian: Some more text
Now we turn our attention to the proof of Theorem~1.4.
We deduce (1.4) from a well-known determinant evaluation 
(see the proof of Lemma~2.2) and Proctor's 
formula (1.1), using the Factorization Theorem for perfect matchings
of \cite{\Ci1}.

\topinsert
\centerline{\mypic{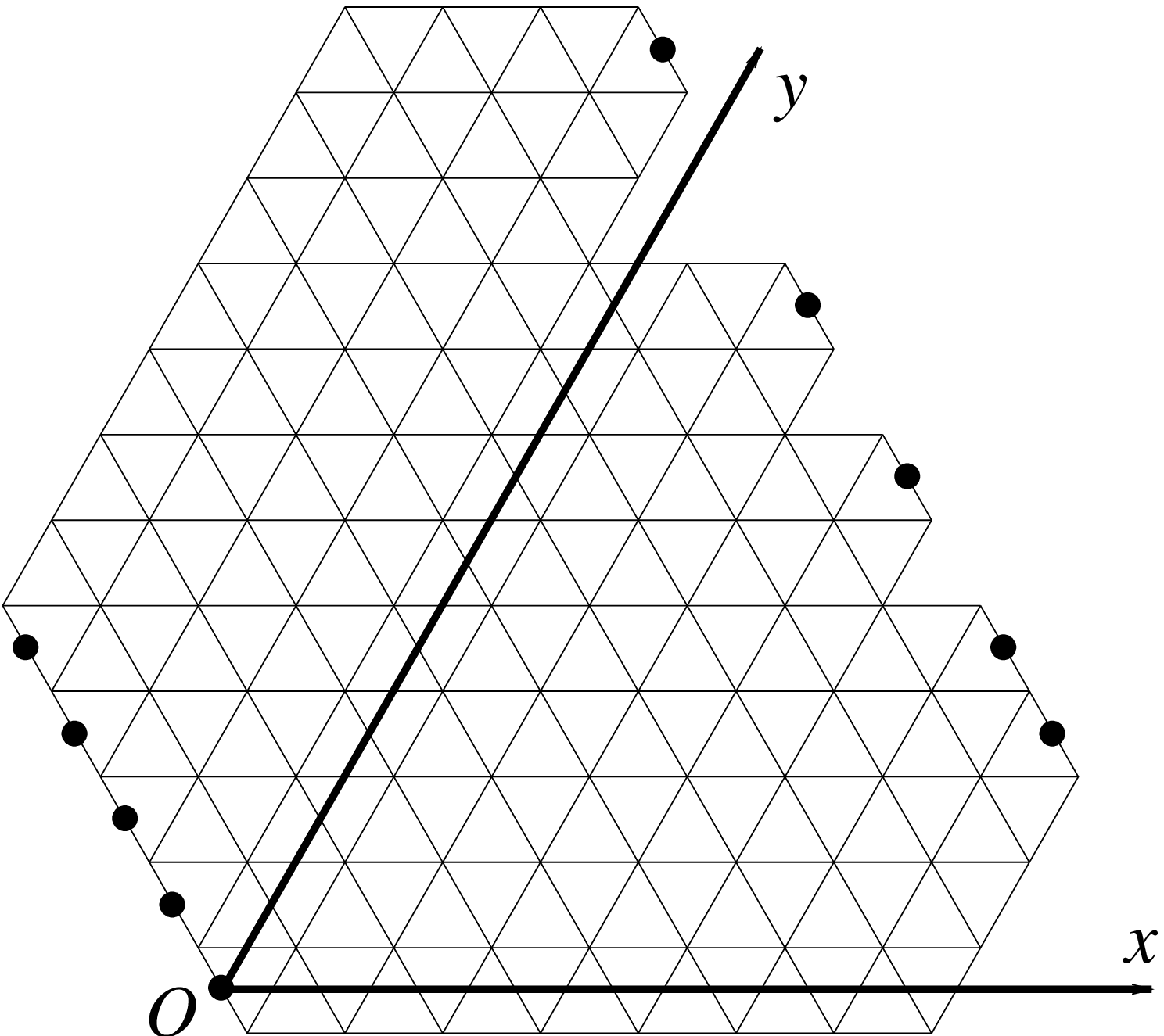}}
\centerline{Figure~2.5. {\rm $R(5,10,(7,6,4,2,-1))$.}}
\endinsert

Let $m$ and $n$ be nonnegative integers, and let ${\bold l}=(l_1,l_2,\dotsc,l_n)$,
$0\leq l_1\leq m$, $l_1>\cdots>l_n$ be a list of integers (some of which may be 
negative). We define the regions $R(n,m,{\bold l})$ as follows. Consider an oblique
coordinate system on the triangular lattice, centered at the midpoint of a lattice
segment facing northeast and having the $x$-axis horizontal and the $y$-axis parallel 
to the lattice
lines going from southwest to northeast (see Figure~2.5). 
Consider, on the one hand, the points $(-i+1,i-1)$,
$i=1,2,\dotsc,n$, and on the other hand the points
$(l_j,m-l_j)$, $j=1,2,\dotsc,n$. We construct $R(n,m,{\bold l})$ so that its
tilings are in bijection with the families of non-intersecting lattice paths with
these starting and ending points. It is easy to see that this determines 
$R(n,m,{\bold l})$ to be the hexagon with side lengths 
$n,l_1,m-l_1,l_1-l_n+1,l_n+n-1,m-l_n-n+1$ (in anticlockwise order, starting with
the southwestern
side), and having triangular dents along the northeastern side at the lattice segments
with midpoints not among $(l_j,m-l_j)$, $j=1,2,\dotsc,n$. ($R(5,10,(7,6,4,2,-1))$
is shown in Figure~2.5.)

\proclaim{Lemma 2.2} 
$$L(R(n,m,{\bold l}))=
\frac {\prod _{1\le i<j\le n} ^{}(l_i-l_j)\prod _{i=1} ^{n}(m+i-1)!} 
{\prod _{i=1} ^{n}(l_i+n-1)!\prod _{i=1} ^{n}(m-l_i)!}.$$
\endproclaim

\pf By construction and by the Lindstr\"om-Gessel-Viennot theorem, we have that 
$$L(R(n,m,{\bold l}))=\det \left(\binom{m}{l_j+i-1}\right)_{1\leq i,j\leq n}.\tag2.3$$
This determinant can be evaluated, e.g., by means of \cite{\KratBN,
(3.12)} with $A=m-1$ and $L_i=l_i-1$, $i=1,2,\dots,n$,
and one obtains the formula in the statement of the Lemma.
%By transposing, the $(i,j)$-entry becomes
%
%$$\binom{m}{l_i+j-1}=\binom{m}{\lambda_i+j-i}=e_{\lambda_i+j-i}(1^m),\tag2.3$$
%
%where $e_k(x_1,x_2,\dotsc)$ is the $k$-th elementary symmetric function, and the
%symbol $1^m$ in the argument stands for the specialization $x_1=\cdots=x_m=1$,
%$x_i=0$, $i>m$. By the N\"agelsbach--Kostka identity (see e.g. \cite{\Macd, I
%\S3, (3.5)}), the
%determinant in (2.3) equals then $s_{\lambda'}(1^m)$, where $s_\mu$ is the Schur
%function indexed by partition $\mu$. Since the length of $\lambda'$ is $l_1\leq m$, 
%by the formula for the principal specialization of Schur functions 
%(see \cite{\Macd, I \S3, Ex. 1}) we obtain the statement of the Lemma.
$\square$

\smallpagebreak
{\it Proof of Theorem 1.4.} 
Consider the region $H_1(a,b,c)$ and weight by $1/2$ the $a$ tile positions required
to be weighted so by $L^*$ (these are marked in Figure~2.6). Draw a line $l$ through
the centers of the marked lozenges, and reflect $H_1(a,b,c)$ across $l$. The union
$U$ of $H_1(a,b,c)$ with its mirror image is precisely the region $R(2c,a,{\bold l})$, 
where ${\bold l}=(b,b-1,\dotsc,b-c+1,a-c,a-c-1,\dotsc,a-2c+1)$. Applying Lemma~2.2
one obtains, after some manipulations, that
$$L(U)=\prod_{i=1}^a\left[\prod_{j=1}^{b-a}\left(\frac{c+i+j-1}{i+j-1}\right)^2
\prod_{j=b-a+1}^b\frac{2c+i+j-1}{i+j-1}\right].\tag2.4$$
The region $U$ is symmetric about $l$, so we can apply to it the Factorization Theorem
of \cite{\Ci1} (see \cite{\CiKrAA} for a phrasing of it in terms of lozenge tilings). 
Following the prescriptions in the statement of this theorem, 
cut $U$ along the zig-zag line following lattice segments
just above $l$ (this is shown as a thick line in Figure~2.6), and denote the pieces
above and below the cut by $U^+$ and $U^-$, respectively. In $U^-$, 
weight the tile 
positions just below the cut by $1/2$. Since $l$ cuts through $2a$
unit triangles, the Factorization Theorem yields
$$L(U)=2^aL(U^+)L^*(U^-).\tag2.5$$
However, $U^-$ is by construction just $H_1(a,b,c)$. Moreover, in $U^+$ there is a row
of forced tiles (shaded in Figure~2.6), and the region left upon their removal is
congruent to $H_1(a,b-1,c)$. Solving for $L^*(U^-)$ in (2.5) and using 
formulas (2.4) and (1.1), one obtains for $L^*(H_1(a,b,c))$ the product expression
(1.4). $\square$

\smallpagebreak
Even though the region $H_{\text a}(a,b,c)$ looks different from the case
$b\leq a\leq c$ of $H_{\text d}(a,b,c)$, it turns out that their tiling enumerations 
amount to evaluating the same determinant.

\proclaim{Lemma 2.3} For $b\ge a+c-1$ we have 
$L(H_{\text a}(a,b,c))=\det A_a(b-a,c)$.
\endproclaim

\topinsert
\twoline{\mypic{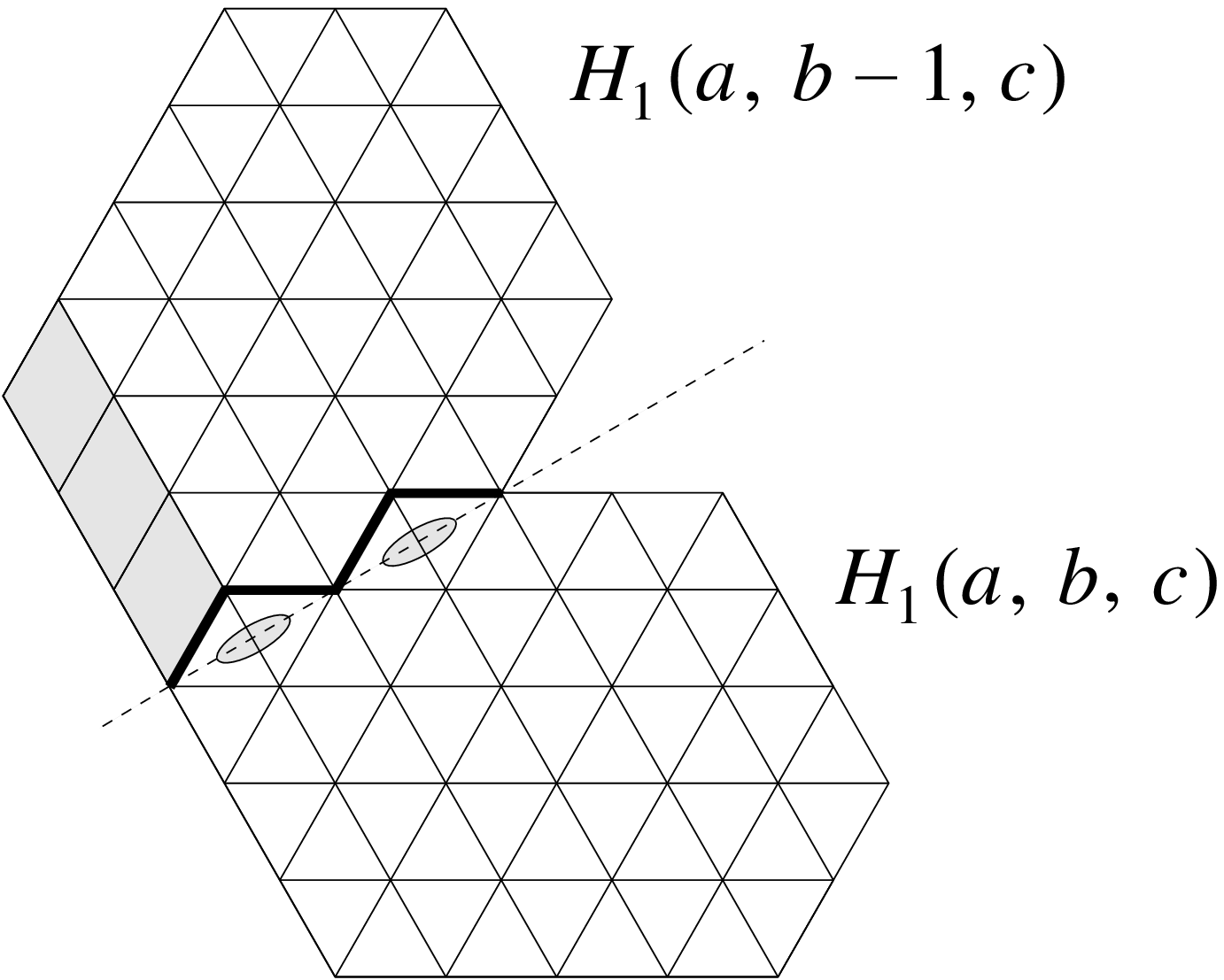}}{\mypic{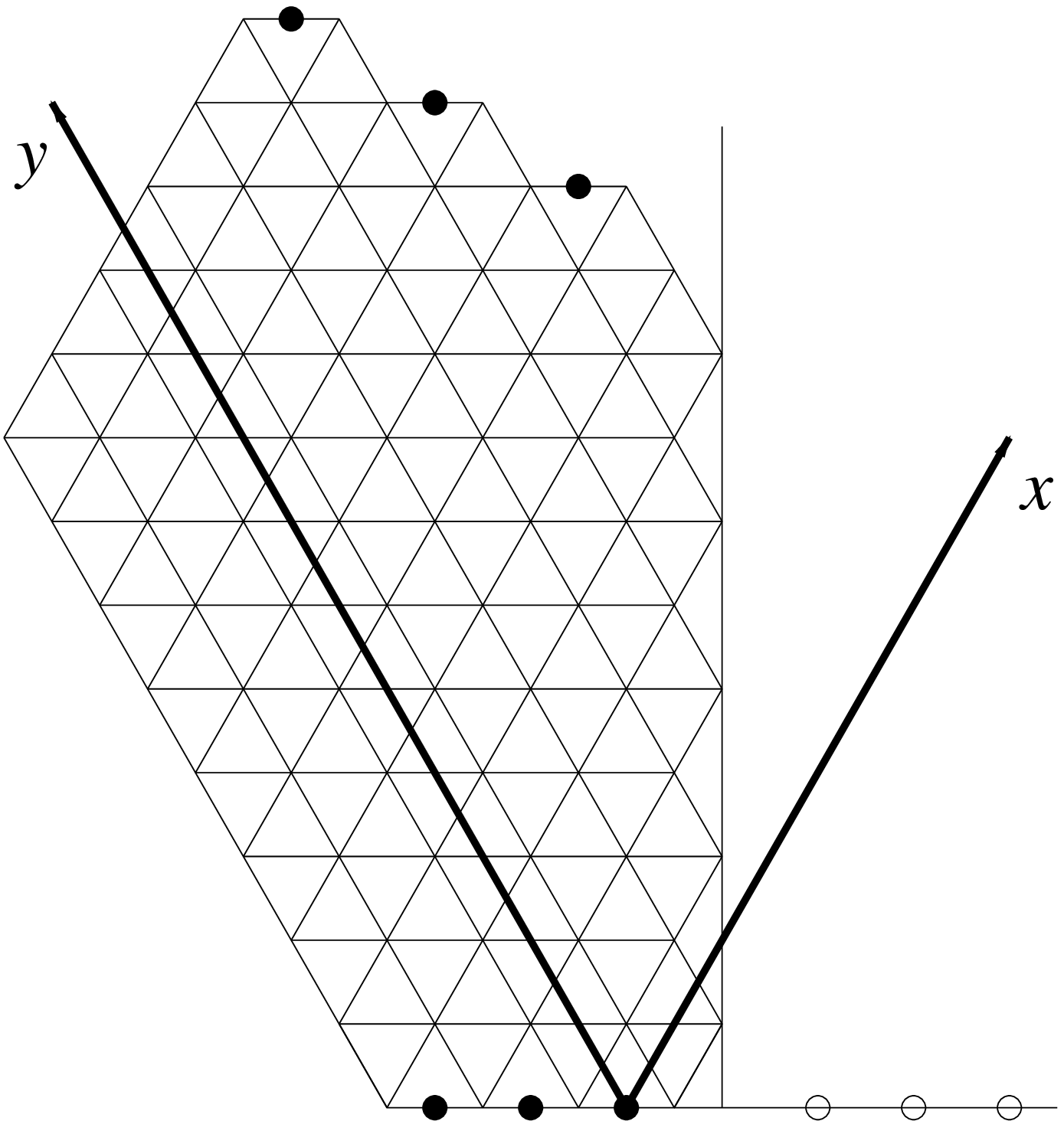}}
\twoline{Figure~2.6.{\rm \ \ \ \ \  }}{\ \ \ \ \ \ \ \ \ \ Figure~2.7.{\rm }}
\endinsert

\pf Rotate the region $H_{\text a}(a,b,c)$ clockwise by 
$60^\circ$, so that it is positioned as in Figure~2.7. 
By the bijection between
tilings and lattice paths, each tiling is identified with a family of 
non-intersecting
lattice paths with starting points $(-i+1,i-1)$, ending points
$(c-j+1,b-a+2j-1)$, $i,j=1,2,\dotsc,a$, and such that all lattice paths stay strictly above
the line $y=x-2$. Just as in the proof of Theorem~1.1, the 
Lindstr\"om-Gessel-Viennot theorem implies that
$$L(H_{\text a}(a,b,c))=
\det \left(\binom{b+c-a+j}{b-a+2j-i}-\binom{b+c-a+j}{b-a+2j+i}
\right)_{1\leq i,j\leq a}.$$
The above determinant is readily recognized 
as $A_a(b-a,c)$. $\square$

\smallpagebreak
{\it Proof of Theorem 1.6.} This follows directly from Lemma 2.3 and Theorem 1.10.
$\square$

\smallpagebreak
{\it Proof of Theorem 1.7.} We use the bijection between tilings and
families of non-inter\-sect\-ing lattice paths from the 
%Christian: Now it is not the preceding proof anymore ...
%preceding 
proof of Lemma~2.3.
In addition, we prepend $(2i-1)$ vertical steps to the $i$-th path.
Thus we obtain families $\Cal P'$ of non-inter\-sect\-ing lattice paths,
with starting points $(-i+1,-i)$, $i=1,2,\dots,a$, ending points
$(c-j+1,b-a+2j-1)$, $j=1,2,\dotsc,a$, 
and such that all lattice paths stay strictly above
the line $y=x-2$. For the rest of the proof one follows the arguments
in the proof of Theorem~1.2, which have to be adjusted only
insignificantly. $\square$

\smallpagebreak
{\it Proof of Proposition 1.8.} Let $T$ be a tiling of $H_{\text o}(a,b,a)$. Consider
the $b$ tiles containing the lattice segments on the bottom part of its boundary.
Because of forcing, there is precisely one dent in the upper boundary of the union of
these tiles (see Figure~2.8). This dent has to be covered by some other tile $t$
(shaded dark in Figure~2.8), which 
in turn forces $b$ more tiles in place. Thus, a subregion congruent to $H(1,b,1)$ at
the bottom of $H_{\text o}(a,b,a)$ ends up being tiled by the restriction of the
tiling $T$. Since there are $b+1$ tilings of $H(1,b,1)$ (corresponding to the $b+1$
possible positions of $t$), this implies
$$L(H_{\text o}(a,b,a))=(b+1)\,L(H_{\text o}(a-1,b,a-1)).$$
Repeated application of this gives the statement of the Proposition. $\square$

\topinsert
\centerline{\mypic{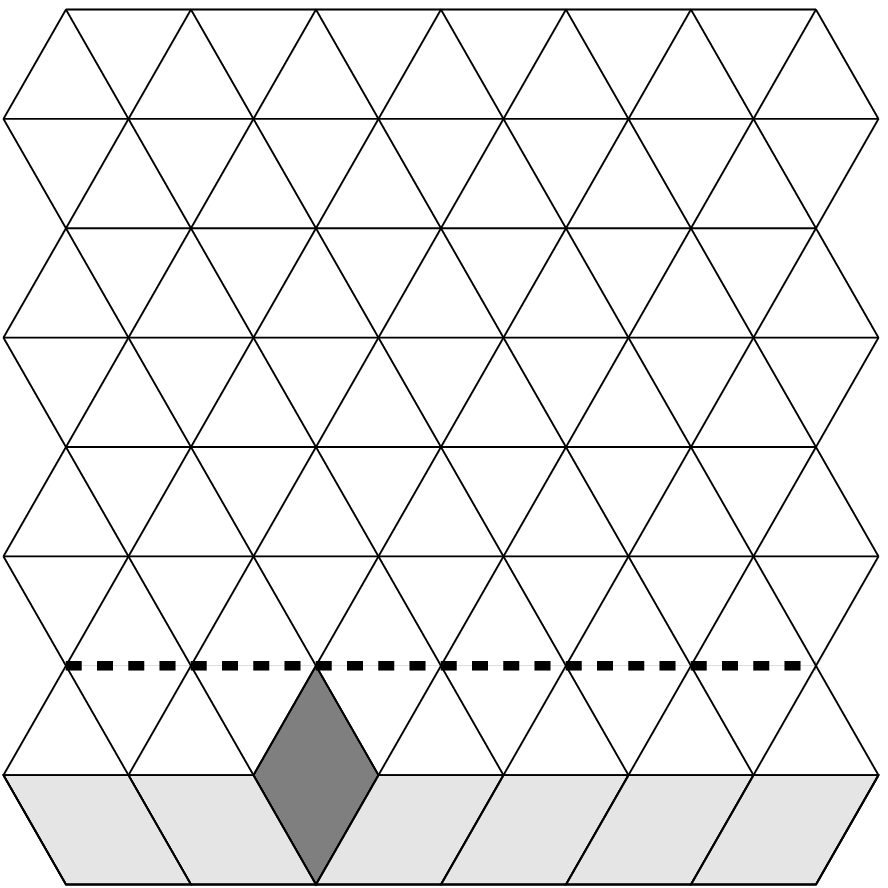}}
\centerline{Figure~2.8. {\rm }}
\endinsert

\mysec{Appendix}

\proclaim{Conjecture A.1}
The number of lozenge tilings of the region $H_{\text n}(x,m+y,x+m-y)$ 
{\rm(}see Figure~{\rm1.7(b)} for an example{\rm)}
is equal to
$$\multline 
\prod _{i=1} ^{m}\frac {(x+i)!} {(x-i+m+y+1)!\,(2i-1)!}
\prod _{i=m+1} ^{m+y}\frac {(x+2m-i+1)!}
{(2m+2y-2i+1)!\,(m+x-y+i-1)!}\\
\times
{2^{\binom {m}2 + \binom y2 }}
      \prod_{i = 1}^{m-1}i!   
      \prod_{i = 1}^{y-1}i!   
      \prod_{i \ge 0}^{}
        ({ \textstyle x+i+{\frac{3}{2}} }) _{m-2i-1}    
      \prod_{i \ge 0}^{}
        ({ \textstyle   x - y+{\frac{5}{2}} + 3 i}) _{ \left \lfloor
         {\frac{3 y}{2}} -{\frac{9 i}{2}}  \right \rfloor-2}\\    
\times      \prod_{i \ge 0}^{}
        ({ \textstyle   x + {\frac{3 m}{2}} - y  + \left \lceil
         {\frac{3 i}{2}} \right \rceil+\frac {3} {2}}) 
      _{ 3 \left \lceil {\frac{y}{2}}
         \right \rceil - \left \lceil
         {\frac{9 i}{2}} \right \rceil -2}    
           \prod_{i \ge 0}
        ^{}
        ({ \textstyle {  x+ {\frac{3 m}{2}} - y + \left \lfloor
         {\frac{3 i}{2}} \right \rfloor+2}}) _{ 3 \left \lfloor
         {\frac{y}{2}} \right \rfloor - \left \lfloor {\frac{9 i}{2}} \right
         \rfloor-1}\\
\times      \prod_{i \ge 0}^{}
        ({ \textstyle  x+m - \left \lfloor {\frac{y}{2}} \right
         \rfloor}+i+1) _{  2 \left \lfloor {\frac{y}{2}} \right
         \rfloor-m - 2 i }
     \prod_{i \ge 0}
        ^{}
        ({ \textstyle  x + \left \lfloor {\frac{y}{2}} \right \rfloor+i+2})
         _{m - 2 \left \lfloor {\frac{y}{2}} \right \rfloor-2i-2} 
\\\times
{\frac{ \dsize
      \prod_{i = 0}^{y}
        ({ \textstyle x - y+3i+1}) _{m + 2 y-4i}    
      \prod_{i = 0}^{ \left \lceil {\frac{y}{2}} \right \rceil-1}
        ({ \textstyle x+m - y+i+1}) _{3 y-m-4i}    
   }
{\dsize
      \prod_{i \ge 0}^{}
        ({ \textstyle x+ {\frac{m}{2}}  - {\frac{y}{2}}+i+1}) _{y-2i}\,
  ({ \textstyle  x + {\frac{m}{2}}-
          {\frac{y}{2}}+i+{\frac{3}{2}}}) _{y-2i-1}   }}\\
\times\frac {\dsize
      \prod_{i = 0}^{y}
        ({ \textstyle x+i+2 }) _{2m - 2 i - 1}    }
 {  ({ \textstyle  x + y+2}) _{ m - y-1}  \,(m+x-y+1)_{m+y} }.
\endmultline\tag A.1$$
Here, shifted factorials occur with positive as well as with negative
indices. The convention with respect to which these have
to be interpreted is
$$(\alpha)_k:=\cases \alpha(\alpha+1)\cdots(\alpha+k-1)&\text {if
}k>0,\\
1&\text {if }k=0,\\
1/(\alpha-1)(\alpha-2)\cdots(\alpha+k)&\text {if }k<0.
\endcases$$ 
All products $\prod
_{i\ge0} ^{}(f(i))_{g(i)}$ in {\rm(A.1)} have to interpreted as the products over
all $i\ge0$ for which $g(i)\ge0$.
\endproclaim

For a proof one could try to proceed as follows:
first we introduce nonintersecting lattice paths, with
starting points along the bottom, and end points along the
northeastern and northwestern zig-zag lines. On
introducing a suitable coordinate system, the starting points can be
represented as $A_i=(-i,i)$, $i=1,2,\dots,m+y$, and the end points as
$E_i=(x-i,2i)$, $i=1,2,\dots,m$, $E_i=(m+y-2i+1,m+x-y+i)$,
$i=m+1,m+2,\dots,m+y$. The corresponding 
Lindstr\"om--Gessel--Viennot determinant is
$$ \det_{1\le i,j\le m+y}\left(\left\{\matrix \binom
{x+i}{x-i+j}&i=1,\dots,m\hfill\\
\binom{x+2m-i+1}{m+y-2i+j+1}&i=m+1,\dots,m+y
\endmatrix\right\}\right).\tag A.2$$
The task is to evaluate this determinant. In principle, after having
taken suitable factors out of the determinant (so that the new
determinant is a polynomial in $x$), the
``identification of factors" method, as described in Section~2.4 of
\cite{\KratBN}, should be capable of accomplishing the determinant
evaluation.

\proclaim{Conjecture A.2}
The weighted count of lozenge tilings of the region 
$H_{\text n}(x,m+y,x+m-y)$, where the lozenges along
the two zig-zag lines are weighted by $1/2$
{\rm(}see Figure~{\rm1.7(b)} for an example; the lozenges that are
weighted by $1/2$ are marked by ellipses{\rm)}, 
is equal to
$$\multline 
\prod _{i=1} ^{m}\frac {(x+i-1)!} {(x-i+m+y+1)!\,(2i-1)!}
\prod _{i=m+1} ^{m+y}\frac {(x+2m-i)!}
{(2m+2y-2i+1)!\,(m+x-y+i-1)!}\\
\times
{2^{\binom {m}2 + \binom y2 }} 
      \prod_{i = 1}^{m-1}i!   
      \prod_{i = 1}^{y-1}i!   
      \prod_{i \ge 0}^{}
        ({ \textstyle x+i+{\frac{3}{2}}}) _{m-2 i -1}    
      \prod_{i \ge 0}^{}
        ({ \textstyle x - y+3i+{\frac{7}{2}}}) _{ \left \lceil
          {\frac{3 y}{2}}-{\frac{9 i}{2}}  \right \rceil-4}\\    
\times           \prod_{i \ge 0}
        ^{}
        ({ \textstyle x+ {\frac{3 m}{2}} - y + \left \lfloor
         {\frac{3 i}{2}} \right \rfloor+\frac {3} {2}}) 
    _{ 3 \left \lceil {\frac{y}{2}}
         \right \rceil-\left \lfloor
         {\frac{9 i}{2}} \right \rfloor -1} 
      \prod_{i \ge 0}^{}
        ({ \textstyle x+ {\frac{3 m}{2}} - y + \left \lceil
         {\frac{3 i}{2}} \right \rceil+1  }) 
        _{ 3 \left \lfloor {\frac{y}{2}} \right
         \rfloor- \left \lceil
         {\frac{9 i}{2}} \right \rceil +1}  \\
\times   \prod_{i \ge 0}
        ^{}
        ({ \textstyle  x+m - \left \lfloor {\frac{y}{2}} \right
         \rfloor+i+1}) _{ 2 \left \lfloor {\frac{y}{2}} \right
         \rfloor-m - 2 i  }     \prod_{i \ge 0}
        ^{}
        ({ \textstyle  x + \left \lfloor {\frac{y}{2}} \right \rfloor+i+2})
         _{ m - 2 \left \lfloor {\frac{y}{2}} \right \rfloor-2i-2} \\
\times
{\frac{\dsize
     ({ \textstyle  x - y+{\frac{1}{2}}}) _{ \left \lfloor
      {\frac{m}{2}} \right \rfloor+2y}  ({ \textstyle x +m- y}) _{y+1}  
      \prod_{i = 0}^{y}
        ({ \textstyle x+i+1}) _{2m - 2 i }    
  }
{({ \textstyle
        x+{\frac{m}{2}} - {\frac{y}{2}}+ {\frac{1}{2}} }) 
       _{\left \lfloor {\frac{ 3y}{2}} \right \rfloor} 
       ({ \textstyle x + {\frac{3 m}{2}}  -
      {\frac{y}{2}}+1}) _{y+1}  ({ \textstyle x + {\frac{m}{2}} +
      {\frac{y}{2}}+1}) _{\left \lceil {\frac{y-2 }{2}} \right \rceil}  
}}\\
\times
\frac {  \dsize    \prod_{i = 0}^{y}
        ({ \textstyle x - y+3i+1}) _{ m + 2 y-4i}    
      \prod_{i = 0}^{ \left \lceil {\frac{y}{2}} \right \rceil-1}
        ({ \textstyle x +m- y+i+1}) _{3 y-m-4i}    } 
{  \dsize (m+x-y)_{m+y+1}\,
    ({ \textstyle x + y + \left \lceil {\frac{m}{2}} \right \rceil}) _{
     \left \lfloor {\frac{m}{2}} \right \rfloor- y + 1}  
     \prod _{i=0} ^{\lceil y/2\rceil-1}(x-y+1+3i)      }\\
\times
\frac {1}
 { \dsize     {\prod}_{i \ge 0}^{}
        ({ \textstyle x+ {\frac{m}{2}} - {\frac{y}{2}}+i+1}) _{y-2i}
  ({ \textstyle x + {\frac{m}{2}}-
          {\frac{y}{2}}+i+{\frac{3}{2}} }) _{y-2i-1}  },
\endmultline \tag A.3$$
with the same conventions as in the previous conjecture.
\endproclaim

For a proof, one could again 
introduce nonintersecting lattice paths, with
starting points and end points as before.
The corresponding Lindstr\"om--Gessel--Viennot determinant is
$$ \det_{1\le i,j\le m+y}\left(\left\{\matrix 
\frac {(x+i-1)!\,(x+j/2)} {(x-i+j)!\,(2i-j)!}&i=1,\dots,m\hfill\\
\frac {(x+2m-i)!\,(3m/2+x-y/2-j/2+1/2)} 
{(m+y-2i+j+1)!\,(m+x-y+i-j)!}&i=m+1,\dots,m+y
\endmatrix\right\}\right).\tag A.4$$
The remarks after (A.2) apply also here.

\mysec{References}
{\openup 1\jot \frenchspacing\raggedbottom
\roster

\myref{\Amdeb}
  T. Amdeberhan, Lewis strikes again!, electronic manuscript dated 1997.
(available at http://www.math.temple.edu/$\sim$tewodros/programs/kradet.html).
%\myref{\AnBuAA}
%  G. E. Andrews and W. H. Burge, Determinant identities, {\it Pacific J. Math.} {\bf
%158} (1993), 1--14.
\myref{\Ci1}
  M. Ciucu, Enumeration of perfect matchings in graphs with reflective symmetry, 
{\it J. Combin\. Theory Ser\.~A} {\bf 77} (1997), 67--97.
\myref{\CiucAI} 
M.    Ciucu, Enumeration of lozenge
tilings of punctured hexagons, {\it J. Combin\. Theory Ser.~A\/}
{\bf 83} (1998), 268--272.
\myref{\Cipp1}
  M. Ciucu, Plane partitions I: a generalization of MacMahon's
formula, 
%Christian: Is there any update?
preprint 
(available at http://www.math.gatech.edu/$\sim$ciucu/pp1.ps).
\myref{\CiKrAA}
  M. Ciucu and C. Krattenthaler, The number of centered lozenge tilings of a symmetric 
hexagon, {\it J. Combin\. Theory Ser\.~A } {\bf 86} (1999), 103-126. 
\myref{\CKpp2}
  M. Ciucu and C. Krattenthaler, Plane partitions II: 5 1/2 symmetry classes, 
{\it Advanced Studies in Pure Mathematics}, 
in: Combinatorial Methods in Representation 
Theory, M. Kashiwara, K. Koike, S. Okada, I. Terada, H. Yamada, eds.,
Advanced Studies in Pure Mathematics, vol.~28, RIMS, Kyoto,
2000, pp.~83--103.
\myref{\CiKrAD}
 M.    Ciucu and C. Krattenthaler, A non-automatic (!) application of 
Gosper's algorithm evaluates a determinant from tiling enumeration,
%Christian:
%preprint, 
Rocky Mountain J. Math\. (to appear), {\tt math/0011047}.
\myref{\DT}
  G. David and C. Tomei, The problem of the calissons, {\it Amer\. Math\. 
Monthly} {\bf 96} (1989), 429--431.
\myref{\EisTAA}
T.    Eisenk\"olbl, Rhombus tilings
of a hexagon with two triangles missing on the symmetry axis,
{\it Electron\. J. Combin\.} {\bf 6} {\rm(1)} (1999) \#R30, 19~pp. 
\myref{\Feller}
  W. Feller, ``An introduction to probability theory and its applications,'' Vol.~I, 
John Wiley \& Sons, 1968.
\myref{\FuKrAC} 
M.    Fulmek and C. Krattenthaler,
The number of rhombus tilings of a symmetric hexagon which contain a
fixed rhombus on the symmetry axis, I,
{\it Ann\. Combin\.} {\bf 2} (1998), 19--40.
\myref{\GV}
  I. M. Gessel and X. Viennot, Binomial determinants, paths, and hook length formulae,
{\it Adv\. in Math\.} {\bf 58} (1985), 300--321.
\myref{\GospAB}
 R. W. Gosper, Decision procedure for indefinite hypergeometric
summation, {\it Proc\. Natl\. Acad\. Sci\. USA} {\bf 75} (1978),
40--42.
\myref{\GrKPAA}
 R. L. Graham, D. E. Knuth and O. Patashnik, ``Concrete Mathematics,'' 
Addison-Wesley, Reading, Massachusetts, 1989.
\myref{\Kast}
  P. W. Kasteleyn, The statistics of dimers on a lattice, I: the number of
dimer arrangements on a quadratic lattice, {\it Physica} {\bf 27} (1961), 
1209--1225.
\myref{\KratBD}
  C. Krattenthaler, Determinant identities and a generalization of the number of 
totally symmetric self-complementary plane partitions, {\it Electron\. J. Combin\.} 
{\bf 4} No.~1 (1997), 
\#R27.
\myref{\KratBN}
  C. Krattenthaler, Advanced determinant calculus, {\it 
S\'eminaire Lotharingien 
Combin\.} {\bf 42} ("The Andrews Festschrift") (1999), paper B42q.
\myref{\Ku}
  G. Kuperberg, Symmetries of plane partitions and the permanent-de\-ter\-mi\-nant
me\-thod, {\it J. Combin\. Theory Ser\.~A} {\bf 68} (1994), 115--151.
\myref{\LindAA}
B.    Lindstr\"om, On the vector
representations of induced matroids, {\it Bull\. London
Math\. Soc\.} {\bf 5} (1973) 85--90.
\myref{\MacM}
  P. A. MacMahon, Memoir on the theory of the partition of numbers --- Part V. Partitions
in two-dimensional space, {\it Phil\. Trans\. R. S.}, 1911, A.
\myref{\PeWZAA}
 M.    Petkov\v sek, H. Wilf and D. Zeilberger, 
``A=B,'' A.K. Peters, Wellesley, 1996.
\myref{\Proc88}
  R. A. Proctor, Odd symplectic groups, {\it Invent\. Math\.} {\bf 92} (1988),
307--332.
\myref{\Robb}
  D. P. Robbins, The story of $1$, $2$, $7$, $42$, $429$, $7436,\dotsc,$ {\it Math.
Intelligencer} {\bf 13} (1991), No.~2, 12--19.
\myref{\Stapp}
  R. P. Stanley, Symmetries of plane partitions, {\it J. Combin\. Theory Ser.~A} 
{\bf 43} (1986), 103--113.
\myref{\Stenlp}
  J. R. Stembridge, Nonintersecting paths, Pfaffians and plane partitions, 
{\it Adv\. in Math\.} {\bf 83} (1990), 96--131.

\endroster\par}

\enddocument